 \documentclass[draft, 12pt, onecolumn]{IEEEtran}                                                         

\IEEEoverridecommandlockouts                              

\title{\LARGE \bf
Weakly Coupled Dynamic Program:\\
 Information and Lagrangian Relaxations
}

\author{Fan~Ye, Helin Zhu, and Enlu~Zhou \\
\small School of Industrial and Systems Engineering, Georgia Institute of Technology, Atlanta, GA, 30332 USA\\ fye8@gatech.edu, hzhu67@gatech.edu, enluzhou@isye.gatech.edu}

\usepackage[all]{xy}                    
\usepackage{amscd}
\usepackage{amsmath}
\usepackage{amssymb}                   
\usepackage{amsfonts}                  
\usepackage[amsmath,thmmarks]{ntheorem}

\usepackage{graphicx,epsfig}                  
\usepackage{mathrsfs}

\usepackage{latexsym}
\usepackage{lscape}
\usepackage{rotating}
\usepackage{graphicx}
\usepackage{mathptmx}
 \usepackage{bbm}

\newtheorem*{thm7}{Theorem 7}
\newtheorem*{thm9}{Theorem 9}

\newtheorem{theorem}{Theorem}

\newtheorem{corollary}[theorem]{Corollary}

\newtheorem{lemma}{Lemma}
\newtheorem{assumption}{Assumption}
  \newtheorem{remark}{Remark}

\begin{document}

\maketitle
\thispagestyle{empty}
\pagestyle{empty}
\pagenumbering{arabic}

\begin{abstract}
``Weakly coupled dynamic program'' describes a broad class of stochastic optimization problems in which multiple controlled stochastic processes  evolve independently but subject to a set of linking constraints imposed on the controls.  One feature of the weakly coupled dynamic program  is that it decouples into  lower-dimensional dynamic programs by dualizing the linking constraint via the Lagrangian relaxation, which also yields a bound on the optimal value of the original dynamic program. Together with the Lagrangian bound, we utilize the information relaxation approach that relaxes the non-anticipative constraint on the controls  to obtain a tighter dual bound. We also investigate other combinations of the relaxations and place the resulting bounds in order. To tackle large-scale problems,  we further propose a  computationally tractable method  based on  information relaxation, and provide insightful interpretation and performance guarantee. We implement our method and demonstrate its use through two numerical examples. 
\end{abstract}



\maketitle

%


Many sequential decision making problems under uncertainty are preferably modelled as  Markov decision processes, which in principle can be solved by dynamic programming. However,  solving large-scale dynamic decision making problems via exact dynamic programming is usually intractable  due to the ``curse of dimensionality'', that is, the computational complexity increases exponentially with the dimension of the state space. To address this, many approximate dynamic programming methods  have been proposed such as \cite{bertsekas:2007, chang:2007a, powell:2011, de2003linear}. In particular, various heuristic  policies are derived based on approximations of the value function.

It is worth noting that the problem-specific approximate value may  be derived based on the relaxations of two intrinsic constraints in general stochastic dynamic programs. One is ``budget constraint'' or the feasibility of the control, which means the decision or control should take values in a feasible region. Another constraint is the ``information constraint'' or non-anticipativity of the control policy, that is, the decision should depend on the information up to the time that the decision is made.  These relaxations may lead to a simpler dynamic optimization problem: the first constraint that exists universally in  mathematical programs  can be tackled by the commonly known  Lagrangian relaxation (see, e.g., \cite{bertsekas1982constrained}), which results in an unconstrained stochastic dynamic program that may be easier to solve; the second constraint can be approached by a recently developed technique --``information relaxation'' (see, e.g, \cite{rogers:2007,brown:2010, brown:2011,Desai:2011,ye2013information}), which relaxes the non-anticipativity constraint  on the controls but impose a penalty for such a violation. Since this approach allows the decision to be made based on the future outcome, it involves  scenario-based dynamic programs, which are deterministic optimization problems and may be less complicated than the original stochastic dynamic program.  There also exist  other relaxation methods. For example, the LP-based approximate dynamic programming (ALP) method  proposed by \cite{schweitzer1985generalized,de2003linear} employs a parameterized class of functions to approximate the optimal value based on the linear programming  formulation of the Bellman optimality equation. These relaxations  provide not only  approximate values  that may be used for developing control policies, but also upper bounds (or lower bounds) on the optimal expected rewards (or expected costs). With such a complementary dual bound on the optimal value,  we can easily evaluate the quality of the  policy available and justify the need of improvement, once its induced expected value    gets close to the dual bound.


In this paper we study the interactions between the aforementioned relaxations  in the  weakly coupled  dynamic program (see, e.g.,  \cite{hawkins2003langrangian,adelman2008relaxations}),   which consists of  multiple subproblems that are independent of each other  except for a set of budget or linking constraints on the controls. This broad class of stochastic optimization problems have many interesting and practical applications including multi-armed and restless bandits, resources allocation, network revenue management, and optimal learning (see, e.g., \cite{jones1972dynamic, gittins1979bandit, whittle1988restless, bertsimas2000restless, gocgun2010lagrangian, talluri1998analysis, topaloglu2009using, bertsimas2007learning, caro2007dynamic, frazier2008knowledge}). However,  the exact solution to the weakly coupled dynamic program quickly becomes  intractable  as the number of subproblems increases; therefore, both approximate solution methods and  dual bounds are needed. In particular, \cite{adelman2008relaxations} compared the Lagrangian relaxation and ALP methods: it is shown that the latter approach provides a tighter  bound on the value function, while the Lagrangian relaxation-based bound is easier to compute. In this paper we further investigate the dual bound  by incorporating the information relaxation approach. We find that in principle the information relaxation approach can be used to improve the quality of bounds derived  from  either  Lagrangian relaxation or ALP method. By combining different types of relaxations, we can obtain various dual bounds with different performance guarantee. As a consequence of the quality versus complexity  trade-off, we may start with one relaxation that requires the least computational cost, and based on the empirical bound performance we may decide how much more we should invest to derive better policies or/and tighter bounds.

We consider the information relaxation in weakly coupled dynamic programs  for both discounted infinite-horizon and finite-horizon cases.  There are a couple of related literatures. \cite{Desai:2011}  explored the  theoretical formulation of the information relaxation bound and compared it with the ALP method in  discounted infinite-horizon Markov decision processes. In contrast, we propose a computational method that can be implemented in the discounted infinite-horizon setting.  \cite{kim2013roubustMAB} studied a robust model of the multi-armed bandit using the information relaxation approach.  \cite{brown2013info}  developed a gradient-based penalty method to compute dual bounds on the revenue in an airline network problem, which is  a case of weakly coupled dynamic programs;  their method can also be used in convex stochastic dynamic programs. In this paper we consider more general weakly coupled dynamic programs with  a class of non-convex linking constraints including the discrete-action case. We propose a computationally tractable approach to generate a weaker bound compared with the exact information relaxation bound, but it is still superior to the Lagrangian relaxation bound; therefore, our computational method  can be viewed as an intermediate relaxation  between the Lagrangian relaxation and the exact information relaxation.

We summarize our contributions as follows:
 \begin{enumerate}
 \item We propose a computational method to obtain an information relaxation  bound on the optimal value of the discounted infinite-horizon weakly coupled dynamic programs, which can also be applied on  discounted infinite-horizon Markov decision processes.
 \item We  analyze various combinations of relaxations and place the resulting bounds in order.
\item We compare the respective sufficient and necessary conditions that  Lagrangian relaxation and information relaxation give  tight bounds on the optimal value. We  provide an example where the Lagrangian bound can be arbitrarily loose, whereas its induced information relaxation bound is always tight.
\item We propose a practical method to compute the information relaxation bound for large-scale  weakly coupled dynamic program and demonstrate its performance guarantee. We also provide the technical conditions such that  the relative gap between the exact and practical information relaxation bounds  vanishes as the number of subproblems goes to infinity.
\end{enumerate}

The rest of this paper is organized as follows. In Section \ref{Section:Formulation}, we review the formulation of the weakly coupled dynamic program and  its decomposition using the Lagrangian relaxation approach.  In Section \ref{Section:Information_Bound}, we present the information relaxation-based dual bounds for the infinite-horizon problem, and compare it to the Lagrangian relaxation and ALP method. In Section \ref{Section:Computational_Method}, we address the computational issue of the information relaxation bound in the large-scale setting. We present our numerical studies in Section \ref{Section:Computation}, and provide the concluding remarks in Section \ref{Section:Conclusion}. We put some of the proofs and the results of the finite-horizon problem in Appendix of the online supplement to this paper.

\section{Formulation of the Weakly Coupled Dynamic Program }\label{Section:Formulation}
In this section, we present the general framework of the weakly coupled stochastic dynamic program and the Lagrangian relaxation approach.
\subsection{Problem  Formulation }
Consider  a collection of $N$ projects or subproblems labeled by $n=1,\cdots,N$. The state of each project or subproblem transits independently according to a homogenous transition law  and yields a reward that is dependent only on the individual state and control. However, at each time period there are  constraints imposed on the controls of these projects, which are referred to as the ``linking constraints'' or ``budget constraints''. The underlying probability space is described by $(\Omega, \mathcal{F}, P)$, where $\Omega$ is the set of possible outcomes or scenarios $\omega$, $\mathcal{F}$ is a $\sigma$-algebra containing  the events in $\Omega$, and $P$ is a probability measure.

We use the following notations to describe the mathematical formulation of the weakly coupled stochastic optimization problem.
\begin{enumerate}
\item Time is indexed by $t=0,1,2\cdots$.
\item $\textbf{x}_t=(x^1_t,\cdots,x^N_t)$ is the joint state of the N projects, and it takes value in the state space $\mathcal{X}=\mathcal{X}^1\times \cdots \times \mathcal{X}^N$.
\item $\textbf{a}_t=(a^1_t,\cdots,a^N_t)$ is the control (or decision variable) that takes value in the control (or action) space $\mathcal{A}=\prod_{n=1}^N\mathcal{A}^n$, where $\mathcal{A}^n$ (resp., $\mathcal{A}$) may depend on $x^n_t$ (resp., $\textbf{x}_t$ ), i.e., $a^n_t\in \mathcal{A}^n(x^n_t)$ (resp., $\textbf{a}_t\in \mathcal{A}(\textbf{x}_t)$).

\item The state of $N$-project  transits in a Markovian fashion; in particular, it evolves as $N$ independent Markov decision processes according to a known homogenous transition law
$$P(\textbf{x}_{t+1}|\textbf{x}_t,\textbf{a}_t)=\prod_{n=1}^{N}P_n(x^{n}_{t+1}| x^n_t, a^n_t),$$
where $\{P_n\}_{n=1}^N$ denotes the controlled transition probability of the individual project. Note that each state $\textbf{x}_{t+1}$ depends on the prior control sequence $\textbf{a}(t)\triangleq(\textbf{a}_0,\textbf{a}_1,\cdots,\textbf{a}_{t})$ and the scenario $\omega$, i.e.,  $\textbf{x}_{0}=\textbf{x}_{0}(\omega)$ and $\textbf{x}_{t+1}=\textbf{x}_{t+1}(\textbf{a}(t),\omega)$ for $t\geq0$, where $\omega$  represents the underlying uncertainty.

\item  At period $t$   the control $\textbf{a}_t$ is chosen by the decision maker subject to a set of $L$ time-invariant linking constraints $\sum_{n=1}^N \textbf{B}^{n}(x^n,a^n)\leq \textbf{b}$, where $\textbf{b}\in \mathbb{R}^L$. Denote the feasible control space associated with the state $\textbf{x}_t$ by
\begin{align}
 \bar{\mathcal{A}}(\textbf{x}_t)=\{\textbf{a}_t\in \mathcal{A}(\textbf{x}_t) :~ \textbf{B}(\textbf{x}_t,\textbf{a}_t)\triangleq\sum_{n=1}^N \textbf{B}^{n}(x^n_t,a^n_t)\leq \textbf{b} \}. \label{linking_constraint}
 \end{align}

\item At period $t$ the $n$-th project or subproblem  incurs a reward of $R^n(x^n_t,a^n_t)$.  The total reward  incurred at time $t$ is of the additive form
    $$R(\textbf{x}_t,\textbf{a}_t)\triangleq\sum_{n=1}^N R^n(x^{n}_{t},a^{n}_{t}).$$

\item  Given a scenario $\omega$, the decision maker chooses a sequence of controls $\textbf{a}=(\textbf{a}_0,\textbf{a}_1,\cdots)$, where each $\textbf{a}_t$ takes value in $\bar{\mathcal{A}}(\textbf{x}_t)$. Such a selection is called a control policy, i.e., $\alpha:\Omega\rightarrow \bar{\mathcal{A}}(\textbf{x}_0)\times\bar{\mathcal{A}}(\textbf{x}_1)\times \cdots$. We denote the set of such control policies as  $\bar{\mathbb{A}}$.

\item The filtration $\mathbb{F}=\{\mathcal{F}_0,\mathcal{F}_1,\mathcal{F}_2,\cdots\}$ describes the evolution of the state information,  where $\mathcal{F}_0\triangleq \sigma\{ \textbf{x}_0\}$ and $\mathcal{F}_t\triangleq \sigma\{ \textbf{x}_0,\cdots,\textbf{x}_t, \textbf{a}_0,\cdots,\textbf{a}_{t-1}\}$ for $t\geq 1$. Since the decision maker determines $\textbf{a}_t$   based only on the information known up to period $t$, each $\textbf{a}_t$ is then $\mathcal{F}_t$-measurable; we call such a control policy $\alpha$ to be non-anticipative and denote the set of  non-anticipative policies by $$\bar{\mathbb{A}}_{\mathbb{F}}=\{\alpha\in \bar{\mathbb{A}}|~\alpha ~\text{is non-anticipative}\}.$$



\item The  expected discounted infinite horizon reward induced by a control policy $\alpha$ is
 \begin{align}
V(\textbf{x}_0; \alpha)\triangleq \mathbb{E}\left[\sum_{t=0}^{\infty}\beta^t R\left( \textbf{x}_t,\textbf{a}_t\right)\bigg|\textbf{x}_0\right], \label{context_value_function1}
\end{align}
 where  $\beta\in(0,1)$ is a discount factor, and $\textbf{a}_t$  is selected by $\alpha$  depending on the scenario $\omega$. The objective of the  decision maker is to maximize  the expected infinite horizon reward   over all non-anticipative policies, given the initial condition $\textbf{x}_0\in \mathcal{X}$:
\begin{align}
V(\textbf{x}_0)=\sup_{\alpha\in \bar{\mathbb{A}}_{\mathbb{F}}}V(\textbf{x}_0; \alpha). \label{context_value_function}
\end{align}

\end{enumerate}

To avoid technical complication, we assume that $\{R^n\}_{n=1}^N$ are uniformly bounded on their respective domain (therefore, $V$ is also bounded), and the supremum  in (\ref{context_value_function}) can be achieved (this is the case, for example, when $\mathcal{X}$ and $\mathcal{A}$ are finite). So $V$ is well-defined for all $\textbf{x}_0\in \mathcal{X}$. Thus, the exact solution to (\ref{context_value_function}) can be obtained by solving the following Bellman optimality equation:
\begin{align}
V(\textbf{x}_0)=\max_{\textbf{a}_0\in \bar{\mathcal{A}}(\textbf{x}_0)}\left\{R(\textbf{x}_0,\textbf{a}_0)+\beta \mathbb{E}\left[V(\textbf{x}_1)|\textbf{x}_0,\textbf{a}_0\right]\right\}. \label{context_value_function_exact}
\end{align}
Denote an optimal  stationary and Markov  policy as $\alpha^*=(\alpha^*_{\delta},\alpha^*_{\delta},\cdots)$, where $\alpha^*_{\delta}:\mathcal{X}\rightarrow \mathcal{A}$  satisfies
$$\alpha^*_{\delta}(\textbf{x}_0)\in\arg \max_{\textbf{a}_0\in \bar{\mathcal{A}}(\textbf{x}_0)}\{R(\textbf{x}_0,\textbf{a}_0)+\beta \mathbb{E}\left[V(\textbf{x}_1)|\textbf{x}_0,\textbf{a}_0\right]\}.$$

The standard value iteration or policy iteration algorithm that can be used to solve (\ref{context_value_function_exact}) quickly becomes intractable as $N$ increases, since the size of its state space is   $|\mathcal{X}|=\prod_{n=1}^N |\mathcal{X}^n|$.


\subsection{Lagrangian Relaxation}\label{Section:Lagrangian_Bound}
In this subsection we  consider  the Lagrangian dual of (\ref{context_value_function_exact}) that relaxes the linking constraints on the controls.  The  motivation of relaxing the linking constraint is to  decompose the original high-dimensional problem to several low-dimensional subproblems.

Denote  by    $\mathbb{A}\triangleq\{\alpha:\Omega \rightarrow \mathcal{A}(\textbf{x}_0)\times\mathcal{A}(\textbf{x}_1)\times \cdots \}$, which contains $\mathbb{\bar{A}}$ as a subset. By  dualizing the linking constraint with the Lagrangian multiplier $\boldsymbol\lambda\in \mathbb{R}^L_{+}$, we define  $J^{\boldsymbol\lambda}(\textbf{x}_0)$ for $\textbf{x}_0\in \mathcal{X}$:
\begin{align}
J^{\boldsymbol\lambda}(\textbf{x}_0)\triangleq \max_{\alpha\in \mathbb{A}_{\mathbb{F}}}J^{\boldsymbol\lambda}(\textbf{x}_0;\alpha) \label{context_Lagrangian1},
\end{align}
where
$$J^{\boldsymbol\lambda}(\textbf{x}_0;\alpha)\triangleq\mathbb{E}\left[\sum_{t=0}^{\infty} \beta^t \left( R( \textbf{x}_t, \textbf{a}_t)+{\boldsymbol\lambda}^{\top}\left[\textbf{b}-\textbf{B}(\textbf{x}_t,\textbf{a}_t)\right] \right)\bigg|\textbf{x}_0\right],$$
and $\mathbb{A}_{\mathbb{F}}\triangleq\{\alpha\in \mathbb{A}|~\alpha ~\text{is non-anticipative}\}.$
%
%
%


 We list some properties  of $J^{\boldsymbol\lambda}$ in Lemma \ref{Proposition:LagrangianInequality}; in particular, $J^{\boldsymbol\lambda}$ is an upper bound on $V$ given any $\boldsymbol\lambda\geq 0$, which will be referred to as the ``Lagrangian bound'' in the following.

%
%

\begin{lemma}[Properties of $J^{\boldsymbol\lambda}$] \label{Proposition:LagrangianInequality}
\begin{enumerate}
\item For any $\boldsymbol\lambda\geq 0$, $J^{\boldsymbol\lambda}(\textbf{x})\geq V(\textbf{x})$ for all $\textbf{x}\in \mathcal{X}$.
\item $J^{\boldsymbol\lambda}(\textbf{x})$ is convex and piecewise linear in $\boldsymbol\lambda\geq0$.
\item For all $\textbf{x}\in \mathcal{X}$, $J^{\boldsymbol\lambda}(\textbf{x})$ can  be written as
\begin{align}
J^{\boldsymbol\lambda}(\textbf{x})= \frac{\boldsymbol\lambda^{\top}\textbf{b}}{1-\beta} \boldsymbol+ \sum_{n=1}^{N}H^{\boldsymbol\lambda,n}(x^n),  \label{context_Lagrangian_equality}
\end{align}
where  $H^{\boldsymbol\lambda,n}(x_0^n)$ is the solution to the following  Bellman optimality equation for each $n=1,\cdots, N$:
\begin{equation}
H^{\boldsymbol\lambda,n}(x_0^n)=\max_{a^n_0\in \mathcal{A}^n(x^n_0)} \big \{R^n( x^n_0,a^n_0)-\boldsymbol\lambda^{\top} \textbf{B}^{n}(x^n_0,a^n_0)+\beta \mathbb{E}\left[H^{\boldsymbol\lambda,n}(x_1^n)|x^n_0\right] \big \}. \label{Subproblem_Bellman}
\end{equation}
\end{enumerate}
\end{lemma}

The proof of these results can be found in Theorem 1 and Theorem 2 of Section 2 in \cite{hawkins2003langrangian}, or Proposition 1 and Proposition 2 in \cite{adelman2008relaxations}.

In the case that $\mathcal{X}$ and $\mathcal{A}$ are finite, we may compute the tightest Lagrangian bound over $\boldsymbol\lambda\geq 0$ via a linear program. To be more specific, suppose  $\{\upsilon(\textbf{x}), \textbf{x}\in \mathcal{X}\}$ is  a probability distribution on $\mathcal{X}$, which can be  viewed as the initial distribution of $\textbf{x}_0$.  Let $\upsilon_n(\cdot)$ denote the marginal distribution of  $\upsilon$ with respect to the $n$-th project, i.e,  $\upsilon_n(x^n_0)=\sum_{\{\textbf{x}=(x^1,\cdots, x^n)\in \mathcal{X}: x^n=x^n_0\}}\upsilon(\textbf{x})$. From (\ref{context_Lagrangian_equality}) we define the  Lagrangian bound based on the initial distribution $\upsilon$ as the weighted sum
\begin{align*}
J^{\boldsymbol\lambda}(\upsilon)\triangleq\sum_{\textbf{x}\in \mathcal{X}}\upsilon(\textbf{x})\cdot J^{\boldsymbol\lambda}(\textbf{x})= \frac{\boldsymbol\lambda^{\top}\textbf{b}}{1-\beta} + \sum_{n=1}^N\sum_{x^n\in  \mathcal{X}_n} \upsilon_n(x^n)H^{\boldsymbol\lambda,n}(x^n).
\end{align*}
The optimal $\boldsymbol\lambda^*=\arg \min_{\{\boldsymbol\lambda\geq 0\}} J^{\boldsymbol\lambda}(\upsilon)$ can be determined by the following linear program (with variables $\boldsymbol\lambda$ and $\{H^{n}(\cdot)\}_{n=1}^N$).
\begin{align}
 \min_{\boldsymbol\lambda, H^{n}(\cdot)}~~ &\frac{\boldsymbol\lambda^{\top}\textbf{b}}{1-\beta} +\sum_{n=1}^N\sum_{x^n\in \mathcal{X}_n} \upsilon_n(x^n)H^n(x^n) \label{Lagragian_LP}\\
s.t.~~ \boldsymbol\lambda &\geq 0, \notag\\
H^n(x_0^n)&\geq R^n( x^n_0,a^n_0)-\boldsymbol\lambda^{\top} \textbf{B}^{n}(x^n_0,a^n_0)+\beta\sum_{x_1^n\in  \mathcal{X}^n}P_n(x^n_1|x^n_0,a^n_0) H^n(x^n_1), \notag\\
  &~~~~~~~~~~~~~~~~~~~~~~~~~~~~~~~~~\text{for all}~(x^n_0,a^n_0) ~\text{with}~ a^n_0\in \mathcal{A}^n(x^n_0). \notag
\end{align}

In the continuous-state or continuous-action case, noting that   $J^{\boldsymbol\lambda}(\upsilon)$  is convex in $\boldsymbol\lambda$ with a fixed probability distribution $\upsilon$, the Lagrangian bound $\boldsymbol J^{\boldsymbol\lambda^*}$ may be solved using the stochastic subgradient method (see, e.g., Section 2.2.1 of \cite{hawkins2003langrangian}). We also review the ALP method to derive an upper bound $H^{LP}$ on $V$ and compare  its bound performance with the Lagrangian bound  in Appendix \ref{Appendix:ALP}.

\section{Information Relaxation-based Dual Bound}\label{Section:Information_Bound}
 Information relaxation-based duality method proposed in  \cite{brown:2010,rogers:2007} can be used to compute a dual bound on the optimal value of finite-horizon stochastic dynamic programs. In this section, we propose a  computational approach based on a randomization idea to extend the information relaxation method to the infinite-horizon setting. This computational method is then used to improve the quality of the Lagrangian bound; in some cases this improvement can be significant. We also analyze the conditions that the two bounds equal the optimal value.

 We will use the following notations. Given $T\in \mathbb{N}$, we denote by
 $\mathcal{A}(T)\triangleq\mathcal{A}(\textbf{x}_0)\times \cdots\times \mathcal{A}(\textbf{x}_T)=\prod_{n=1}^N\mathcal{A}^n(T)$. Respectively, we define  $\bar{\mathcal{A}}(T)\triangleq\bar{\mathcal{A}}(\textbf{x}_0)\times \cdots\times \bar{\mathcal{A}}(\textbf{x}_{T}).$

\subsection{Information Relaxation-based Bounds for Discounted Infinite-Horizon Problem}
The Lagrangian relaxation approach in Section \ref{Section:Lagrangian_Bound} relaxes the  feasible set of the controls, where the  term $\sum_{t=0}^{\infty} \beta^t \boldsymbol\lambda^{\top}\big([\textbf{b}-\textbf{B}(\textbf{x}_t,\textbf{a}_t)]\big)$ plays the role of a penalty when the decision   takes value outside the feasible region.   As an alternative relaxation  technique, the ``information relaxation'' relaxes the non-anticipativity constraint on the control policy and impose a class of penalties that penalize this violation.

To begin with, we define the partial sum $M_k$ of a sequence of random variables $\{m_t\}_{t=0}^\infty$ with discount factor $\beta\in(0,1)$, that is,
 \begin{align}
M_k(\textbf{a},\omega)\triangleq \sum_{t=0}^{k} \beta^{t+1} m_t(\textbf{a}(t),\omega) , ~~k=0,1,\cdots, \label{martingale_difference_original}
\end{align}
where $m_t$ depends on the scenario $\omega$  and the decisions up to time $t$, i.e., $\textbf{a}(t)=(\textbf{a}_0,\cdots,\textbf{a}_t)$. In particular, we consider a special form of $m_t$:
\begin{align*}
 m_t(\textbf{a}(t),\omega)=\Delta H(\textbf{x}_{t+1}, \textbf{x}_t,\textbf{a}_t) \triangleq H( \textbf{x}_{t+1})-  \mathbb{E}[ H(\textbf{x}_{t+1})|\textbf{x}_t, \textbf{a}_t],~ \textbf{a}_t\in \mathcal{A}(\textbf{x}_t), ~t=0,1,\cdots, 
 \end{align*}
where  $H\in\mathcal{D}\triangleq\{ H:\mathcal{X} \rightarrow \mathbb{R}|H~\text{is  measurable and bounded}\}$. Note that with a control policy $\alpha\in \mathbb{\bar{A}}_\mathbb{F}$,  $\{ M_k(\alpha,\omega)\}_{k=0}^{\infty}$ is an $\mathbb{F}$-martingale, since $\{m_t\}_{t=0}^{\infty}$ is an $\mathbb{F}$-martingale difference sequence. In particular, $\mathbb{E}[M_k(\alpha,\omega)|\textbf{x}_0]=0$ for any $\alpha\in \mathbb{\bar{A}}_\mathbb{F}$.

We then  consider the discounted infinite sum of $m_t$, that is,
$$M(\textbf{a},\omega)\triangleq \sum_{t=0}^{\infty}\beta^{t+1} m_t(\textbf{a}(t),\omega).$$

We can show that $M(\textbf{a},\omega)$   is well defined for any $\textbf{a}$ and $\omega$ given $H\in \mathcal{D}$  , i.e., $|H(\cdot)|<\Lambda$ for some $\Lambda>0$; the sequence $\{M_k\}_{k=0}^\infty$ is then uniformly bounded  for all $k\geq 0$, since
$$ |M_k(\textbf{a},\omega)| \leq \sum_{t=0}^{k} \beta^{t+1}  |\Delta H(\textbf{x}_{t+1}, \textbf{x}_t,\textbf{a}_t)|\leq \frac{2  \Lambda }{1-\beta}~~\text{for all}~\omega\in \Omega~\text{and}~\textbf{a}_t\in \mathcal{A}(\textbf{x}_t), ~t=0,1,\cdots,k.$$
Therefore,  $M(\textbf{a},\omega)\triangleq\lim_{k\rightarrow \infty}M_k(\textbf{a},\omega)$ is well-defined for every $\textbf{a}$ and $\omega$. In particular, $\mathbb{E}[M(\alpha,\omega)|\textbf{x}_0]=\lim_{k\rightarrow \infty}\mathbb{E}[M_k(\alpha,\omega)|\textbf{x}_0]=0$ for $\alpha\in \mathbb{\bar{A}}_\mathbb{F}$ due to the dominated convergence theorem.

Suppose now $\textbf{a}$ is the control sequence selected by a  policy  $\alpha\in \mathbb{\bar{A}}_\mathbb{F}$, and let  $M(\alpha,\omega)= \sum_{t=0}^{\infty}\beta^{t+1} \Delta H(\textbf{x}_{t+1}, \textbf{x}_t,\textbf{a}_t)$ with $H\in \mathcal{D}$. Then

\begin{align}
V(\textbf{x}_0;\alpha)=&\mathbb{E}\left[\sum_{t=0}^{\infty}\beta^t R(\textbf{x}_t,\textbf{a}_t)\bigg|\textbf{x}_0\right] -\mathbb{E}\left[M(\alpha,\omega)|\textbf{x}_0\right]\notag \\
=&\mathbb{E}_0\left[\sum_{t=0}^{\infty}\beta^t \big(R( \textbf{x}_t,\textbf{a}_t)-  \beta (H( \textbf{x}_{t+1})- \mathbb{E}[ H(\textbf{x}_{t+1})|\textbf{x}_t, \textbf{a}_t]) \big)\right] \notag\\
=&H( \textbf{x}_{0})+\mathbb{E}_0\left[ \sum_{t=0}^{\infty}\beta^t \big(R( \textbf{x}_t,\textbf{a}_t)+ \beta \mathbb{E}[ H(\textbf{x}_{t+1})|\textbf{x}_t, \textbf{a}_t]- H( \textbf{x}_{t})\big)\right], \label{weak_duality_equality1}
\end{align}
where $\mathbb{E}_0\big[\cdot]=\mathbb{E}\big[\cdot|\textbf{x}_0]$. The first equality  holds due to $\mathbb{E}[M(\alpha,\omega)|\textbf{x}_0]$=0 for $\alpha\in \mathbb{A}_{\mathbb{F}}$,  the second equality holds due to the definition of $\Delta H$, and the last equality holds since  $\sum_{t=0}^{\infty} \beta^{t}H( \textbf{x}_{t})$ and  $\sum_{t=1}^{\infty}\beta^{t+1}\mathbb{E}[ H(\textbf{x}_{t+1})|\textbf{x}_t, \textbf{a}_t]$ are absolutely convergent for all $\omega\in\Omega$ and $\textbf{a}\in \mathbb{A}$. 

To develop a computational method that reduces the infinite sum inside the conditional expectation in (\ref{weak_duality_equality1}) to a finite sum, we consider a random time $\tau$ (see, .e.g., \cite{fox1989simulating}) that  is independent of $\{\mathcal{F}_t, ~t=0,1,\cdots\}$, and $\tau$ is of  geometric distribution with parameter $\beta$, i.e.,
\begin{align}
P(\tau=t)= (1-\beta)\beta^t, ~~t=0,1,\cdots. \label{context_random_time}
\end{align}
A complete definition of $\tau$  is in Appendix B.
Noting  that
$P(t\leq \tau)=\mathbb{E}\left[\mathbbm{1}_{\{t\leq\tau\}}\right]=\beta^t$,  we can rewrite the second term in (\ref{weak_duality_equality1}) as
\begin{align*}
&\mathbb{E}_0\left[\sum_{t=0}^{\infty}\mathbb{E}\left[\mathbbm{1}_{\{t\leq\tau\}}\right] \cdot \left(R( \textbf{x}_t,\textbf{a}_t)+\beta \mathbb{E}[ H(\textbf{x}_{t+1})|\textbf{x}_t, \textbf{a}_t]- H(\textbf{x}_{t})\right)\right]\\
  =&\mathbb{E}_0\left[\sum_{t=0}^{\infty}\mathbbm{1}_{\{t\leq\tau\}} \cdot \left(R( \textbf{x}_t,\textbf{a}_t)+\beta \mathbb{E}[ H(\textbf{x}_{t+1})|\textbf{x}_t, \textbf{a}_t]- H(\textbf{x}_{t})\right)\right]\\
 =& \mathbb{E}_0\left[\sum_{t=0}^{\tau} \left(R( \textbf{x}_t,\textbf{a}_t)+\beta \mathbb{E}[ H(\textbf{x}_{t+1})|\textbf{x}_t, \textbf{a}_t]- H(\textbf{x}_{t})\right) \right],
\end{align*}
where the first equality holds due to the Fubini's theorem, noting that the boundedness of $R$ and $H$ implies the integrability of the integrand in $\mathbb{E}_0[\cdot]$.

Based on this transformation, we introduce an operator $\mathcal{L}:\mathcal{D}\rightarrow \mathcal{D}$
\begin{align}
\mathcal{L}H(\textbf{x}_0)\triangleq H(\textbf{x}_{0})+\mathbb{E}\left[\max_{\textbf{a}\in\mathcal{\bar{A}}(\tau) }\{I_H(\textbf{a},\omega,\tau)\}\bigg|\textbf{x}_0\right], \label{context_dual_rep}
\end{align}
where
\begin{equation}
I_{H}(\textbf{a},\omega,\tau)\triangleq \sum_{t=0}^{\tau}\big(\textbf{R}(\textbf{x}_t,\textbf{a}_t)+\beta \mathbb{E}[ H(\textbf{x}_{t+1})|\textbf{x}_t,\textbf{a}_t]- H(\textbf{x}_{t})\big). \label{context_inneropt}
\end{equation}
Note that given $\omega\in\Omega$ and $\tau\in \mathbb{N}$, the dependence of $I_{H}$ on $\textbf{a}$ is only through the first $\tau+1$ actions, namely, $\textbf{a}(\tau)$. Thus, $\max_{\textbf{a}\in\mathcal{\bar{A}}(\tau) }\{I_H(\textbf{a},\omega,\tau)\}$ is short for $\max_{\textbf{a}(\tau)\in\mathcal{\bar{A}}(\tau) }\{I_H(\textbf{a}(\tau),\omega,\tau)\}$, which is referred to as the \emph{inner optimization problem}.   The conditional expectation in (\ref{context_dual_rep}) is now taken with respect to both the random outcome $\omega$ and the random time $\tau$. We can better interpret this conditional expectation via Monte Carlo simulation: in each trial of simulation we first generate a sample of random horizon $\tau$ (that is finite) and a scenario $\omega$, i.e., the underlying uncertainty that affects the evolution of $\{\textbf{x}_t\}_{t=0}^{\tau}$; based on these realizations we  maximize $I_{H}(\textbf{a},\omega,\tau)$ subject to $\textbf{a}\in \mathcal{\bar{A}}(\tau)$ and the state evolution $\{\textbf{x}_t(\textbf{a}(t-1),\omega)\}_{t=0}^{\tau}$. We show that the estimator $\max_{\textbf{a}\in\mathcal{\bar{A}}(\tau)}\{I_H(\textbf{a},\omega,\tau)\}$ has finite  mean and variance in Appendix B.



%

 We next show for any $H\in\mathcal{D}$,  the optimal value $V$ is upper bounded by $\mathcal{L}H$, which will be referred to as the ``information relaxation bound''.  The relaxed information is reflected in the scenario-based inner optimization problem, while   $M(\textbf{a},\omega)= \sum_{t=0}^{\infty}\beta^{t+1} \Delta H(\textbf{x}_{t+1}, \textbf{x}_t,\textbf{a}_t) $  induced by the function $H$ plays the role of a penalty: if $H$ is chosen to be $V$, then the upper bound $\mathcal{L}H$ is  tight, i.e, $\mathcal{L}H=V$.

\begin{theorem}[Information Relaxation Bound] \label{Theorem:context_duality}
Let $\tau$ be a random time of  geometric distribution with parameter $\beta$ and it is independent of $\{\mathcal{F}_t, ~t=0,1,\cdots\}$.  Then
\begin{enumerate}
\item[(a)] (Weak Duality) For any $H\in \mathcal{D}$, $ V(\textbf{x})\leq \mathcal{L}H(\textbf{x})$ for all $\textbf{x}\in \mathcal{X}$.

\item[(b)] (Tighter Bound) For any $H\in\mathcal{D}^\ast$, where
$$\mathcal{D}^\ast\triangleq\{H\in \mathcal{D}: R(\textbf{x}_0,\textbf{a}_0)+\beta \mathbb{E}[ H(\textbf{x}_{1})|\textbf{x}_0,\textbf{a}_0]\leq H(\textbf{x}_{0}) ~~\text{for all}~ \textbf{x}_0\in  \mathcal{X}  ~\text{and}~ \textbf{a}_0\in \bar{\mathcal{A}}(\textbf{x}_0) \},$$
then   $\max_{\textbf{a}\in \mathcal{\bar{A}}(\tau)}\{I_H(\textbf{a},\omega,\tau)\}\leq0$ for every $\omega\in \Omega$ and $\tau\in\mathbb{N}$; consequently,  $\mathcal{L}H(\textbf{x})\leq H(\textbf{x})$   for all $\textbf{x}\in \mathcal{X}$.

\item[(c)]  (Strong Duality)  $V(\textbf{x})=\mathcal{L}V(\textbf{x})$ for all $\textbf{x}\in  \mathcal{X}$.
\end{enumerate}
\end{theorem}
\begin{IEEEproof}
\begin{enumerate}
\item[(a)]  For $\textbf{x}_0\in \mathcal{X}_0$  and $\alpha\in \mathbb{\bar{A}}_{\mathbb{F}}$,
\begin{align*}
V(\textbf{x}_0;\alpha)&= H(\textbf{x}_0)+\mathbb{E}_0\left[I_{H}(\textbf{a},\omega,\tau) \right]\leq H(\textbf{x}_0)+\mathbb{E}_0\left[\max_{\textbf{a}'\in \mathcal{\bar{A}}(\tau)}\left\{I_{H}(\textbf{a}',\omega,\tau)\right\} \right].
\end{align*}
where $\textbf{a}$ is the control sequence selected by $\alpha$. By maximizing $V(\textbf{x}_0;\alpha)$ over $\alpha\in \mathbb{\bar{A}}_{\mathbb{F}}$, the weak duality $V(\textbf{x}_0)\leq \mathcal{L}H(\textbf{x}_0)$ holds.

\item[(b)]  Note that given  any $H\in \mathcal{D}^*$ and $\textbf{x}_t\in \mathcal{X}$,
$R( \textbf{x}_t,\textbf{a}_t)+\beta \mathbb{E}[ H(\textbf{x}_{t+1})|\textbf{x}_t,\textbf{a}_t]- H(\textbf{x}_{t})\leq 0$
for all $\textbf{a}_t\in \bar{\mathcal{A}}(\textbf{x}_t)$. It is straightforward to see that for any $\tau\in \mathbb{N}$ and $\omega\in \Omega$,
$$I_H(\textbf{a},\omega,\tau)=\sum_{t=0}^{\tau}\big(R( \textbf{x}_t,\textbf{a}_t)+\beta \mathbb{E}[ H(\textbf{x}_{t+1})|\textbf{x}_t,\textbf{a}_t]- H(\textbf{x}_{t})\big)\leq 0$$
 for  any $\textbf{a}_t\in \bar{\mathcal{A}}(\textbf{x}_t)$, $t=0,1,\cdots,\tau$. Therefore, for all $\textbf{x}_0\in \mathcal{X}$ we have
 \begin{align*}
 \mathcal{L}H(\textbf{x}_0)=&H(\textbf{x}_{0})+\mathbb{E}\left[\max_{\textbf{a}\in \mathcal{\bar{A}}(\tau)}\left\{I_H(\textbf{a},\omega,\tau)\right\}\bigg|\textbf{x}_0\right] \leq H(\textbf{x}_{0}).
 \end{align*}
Together with the weak duality, we have shown that $V(\textbf{x}_0)\leq\mathcal{L}H(\textbf{x}_0)\leq H(\textbf{x}_0)$.

\item[(c)] Since $V\in \mathcal{D}^*$, the strong duality follows from the last proof by choosing $H=V$.
\end{enumerate}
\end{IEEEproof}

  The function $H\in\mathcal{D}^*$ is sometimes referred to as a ``supersolution'' to the problem (\ref{context_value_function}), and it is a standard result that the optimal value $V$ is upper bounded by a supersolution $H$ (see, e.g., \cite{bertsekas:2007}).  Theorem \ref{Theorem:context_duality}(b) indicates that the scenario-dependent inner optimization problem   of an arbitrary time horizon $\tau$ is upper bounded by zero  provided $H\in \mathcal{D}^*$; therefore, $\mathcal{L}H$ improves the quality of the supersolution $H$ as an upper bound on $V$. The strong duality  implies that we may obtain a tight dual bound, given some approximate function of $V$ that induces a good approximation of $\sum_{t=0}^{\infty}\beta^{t+1} \Delta V(\textbf{x}_{t+1}, \textbf{x}_t,\textbf{a}_t)$.  In addition, Theorem \ref{Theorem:context_duality} is true not only for weakly coupled dynamic program, but also for general discounted infinite-horizon stochastic dynamic program due to the applicable randomization technique.

As a corollary of Theorem \ref{Theorem:context_duality}, we provide the information relaxation-based dual representation of the Lagrangian bound $J^{\boldsymbol\lambda}(\textbf{x})$. Fix $\boldsymbol\lambda\geq 0$ and define the operator $\mathcal{L}^{\boldsymbol\lambda}:\mathcal{D}\rightarrow \mathcal{D}$
{\small\begin{align}
\mathcal{L}^{\boldsymbol\lambda}H(\textbf{x}_0) \triangleq H(\textbf{x}_{0})+\mathbb{E}_0\left[\max_{\textbf{a}\in\mathcal{A}(\tau)} \left\{\sum_{t=0}^{\tau}\left(\textbf{R}(\textbf{x}_t,\textbf{a}_t)+ \boldsymbol\lambda^{\top}[\textbf{b}-\textbf{B}(\textbf{x}_t,\textbf{a}_t)]+\beta \mathbb{E}[ H(\textbf{x}_{t+1})|\textbf{x}_t,\textbf{a}_t]- H(\textbf{x}_{t})\right)\right\}\right]. \label{Lagragian_duality}
\end{align}}
\begin{corollary}\label{Corollary:Dual_Lagrangian}
Suppose $\boldsymbol\lambda\geq 0$. Then
\item[(a)]  For any $H\in \mathcal{D}$, $J^{\boldsymbol\lambda}(\textbf{x})\leq  \mathcal{L}^{\boldsymbol\lambda}H(\textbf{x})$ for all $\textbf{x}\in \mathcal{X}$.
\item[(b)]  $J^{\boldsymbol\lambda}(\textbf{x})=\mathcal{L}^{\boldsymbol\lambda}J^{\boldsymbol\lambda}(\textbf{x})$ for all $\textbf{x}\in  \mathcal{X}$.
\end{corollary}
\begin{IEEEproof}
 Note that the definition of $\mathcal{L}^{\boldsymbol\lambda}H$ parallels that of $\mathcal{L}H$ except for  the one-period reward is
$\textbf{R}( \textbf{x}_t,\textbf{a}_t)+ \boldsymbol\lambda^{\top}[\textbf{b}-\textbf{B}(\textbf{x}_t,\textbf{a}_t)]$ (instead of $\textbf{R}( \textbf{x}_t,\textbf{a}_t)$), and the constraint of the inner optimization problem  is $\mathcal{A}(\tau)$ (instead of $\mathcal{\bar{A}}(\tau)$).  One can directly verify the weak duality, i.e.,  $J^{\boldsymbol\lambda}(\textbf{x})\leq  \mathcal{L}^{\boldsymbol\lambda}H(\textbf{x})$. The strong duality $J^{\boldsymbol\lambda}(\textbf{x})=\mathcal{L}^{\boldsymbol\lambda}J^{\boldsymbol\lambda}(\textbf{x})$ follows from the fact that for   every $\omega\in \Omega$ and $\tau\in\mathbb{N}$,
\begin{align}
\max_{\textbf{a}\in \mathcal{A}(\tau)} \left\{\sum_{t=0}^{\tau}\left(R(\textbf{x}_t,\textbf{a}_t)+\boldsymbol\lambda^{\top}[\textbf{b}-\textbf{B}(\textbf{x}_t,\textbf{a}_t)]+\beta \mathbb{E}[ J^{\boldsymbol\lambda}(\textbf{x}_{t+1})|\textbf{x}_t,\textbf{a}_t]- J^{\boldsymbol\lambda}(\textbf{x}_{t})\right)\right\}=0. \label{inner_LagrangianRelaxation}
\end{align}
\end{IEEEproof}

\subsection{Comparison with the Lagrangian Relaxation}\label{SubSection:ComparingBounds}
In weakly coupled stochastic dynamic program,  the Lagrangian bound $J^{\boldsymbol\lambda}(\textbf{x})=\frac{\boldsymbol\lambda^{\top}\textbf{b}}{1-\beta} + \sum_{n=1}^{N}H^{\boldsymbol\lambda,n}(x^n)$ and the upper bound derived from the ALP method, i.e., $H^{LP}(\textbf{x})=\theta^*+ \sum_{n=1}^{N}H^{LP,n}(x^n)$ (see the definition of $H^{LP}$ in Appendix \ref{Appendix:ALP}) are natural candidates as approximate value functions. It can be shown that  the information relaxation approach indeed improves the performance of both  bounds.

\begin{theorem}\label{Propostion:Lagrangian_Information_Bound}
\begin{enumerate}
\item[(a)] For any $\boldsymbol\lambda\geq0$,  $\mathcal{L}J^{\boldsymbol\lambda}(\textbf{x})\leq J^{\boldsymbol\lambda}(\textbf{x})$ for all  $\textbf{x}\in \mathcal{X}$.

\item[(b)] $\mathcal{L}H^{LP}(\textbf{x})\leq H^{LP}(\textbf{x})$ for all  $\textbf{x}\in \mathcal{X}$.

\item[(c)] If  $H(\textbf{x}_{0})- (R(\textbf{x}_0,\textbf{a}_0)+\beta \mathbb{E}[ H(\textbf{x}_{1})|\textbf{x}_0,\textbf{a}_0] ) \geq  \varepsilon$ for all  $\textbf{x}_0\in \mathcal{X}$  and $\textbf{a}_0\in \bar{\mathcal{A}}(\textbf{x}_0)$, then $\mathcal{L}H(\textbf{x}) \leq H(\textbf{x})- \frac{\varepsilon}{1-\beta}$ for all  $\textbf{x}\in \mathcal{X}$.
\end{enumerate}

\end{theorem}

\begin{IEEEproof}
\begin{enumerate}
\item[(a)]
This is an immediate corollary of Theorem \ref{Theorem:context_duality}(c) since $J^{\boldsymbol\lambda}\in \mathcal{D}^*$ (see Lemma 4(b) in Appendix \ref{Appendix:ALP}). Here we consider an alternative proof based on Corollary \ref{Corollary:Dual_Lagrangian} by showing $\mathcal{L}J^{\boldsymbol\lambda}(\textbf{x})\leq \mathcal{L}^{\boldsymbol\lambda} J^{\boldsymbol\lambda}(\textbf{x})$. Note that for each scenario $\omega$ and $\tau\in \mathbb{N}$,
\begin{align}
0=&\max_{\textbf{a}\in \mathcal{A}(\tau)} \left\{\sum_{t=0}^{\tau}\left(R(\textbf{x}_t,\textbf{a}_t)+\boldsymbol\lambda^{\top}[\textbf{b}-\textbf{B}(\textbf{x}_t,\textbf{a}_t)]+\beta \mathbb{E}[ J^{\boldsymbol\lambda}(\textbf{x}_{t+1})|\textbf{x}_t,\textbf{a}_t]- J^{\boldsymbol\lambda}(\textbf{x}_{t})\right)\right\} \notag\\
\geq&\max_{\textbf{a}\in \mathcal{\bar{A}}(\tau)} \left\{\sum_{t=0}^{\tau}\big(R(\textbf{x}_t,\textbf{a}_t)+\boldsymbol\lambda^{\top}[\textbf{b}-\textbf{B}(\textbf{x}_t,\textbf{a}_t)]+\beta \mathbb{E}[ J^{\boldsymbol\lambda}(\textbf{x}_{t+1})|\textbf{x}_t,\textbf{a}_t]- J^{\boldsymbol\lambda}(\textbf{x}_{t})\big)\right\} \notag\\
\geq&\max_{\textbf{a}\in \mathcal{\bar{A}}(\tau)}\left\{\sum_{t=0}^{\tau}\big(R(\textbf{x}_t,\textbf{a}_t)+\beta \mathbb{E}[ J^{\boldsymbol\lambda}(\textbf{x}_{t+1})|\textbf{x}_t,\textbf{a}_t]- J^{\boldsymbol\lambda}(\textbf{x}_{t})\big)\right\}, \label{WeaklyCoupled_pathwise_inequality}
\end{align}
where the equality follows (\ref{inner_LagrangianRelaxation}), the first inequality holds because $\mathcal{A}(\tau)\supset \mathcal{\bar{A}}(\tau)$, and the second inequality holds since $\boldsymbol \lambda\geq0$ and each $\textbf{b}-\textbf{B}(\textbf{x}_t,\textbf{a}_t)\geq 0$ for $\textbf{a}_t\in \bar{\mathcal{A}}(\textbf{x}_t)$. Hence, $\mathcal{L}J^{\boldsymbol\lambda}(\textbf{x})\leq \mathcal{L}^{\boldsymbol\lambda}J^{\boldsymbol\lambda}(\textbf{x})=J^{\boldsymbol\lambda}(\textbf{x})$ for all $\textbf{x}\in\mathcal{X}$.

\item[(b)]  Note that $H^{LP}\in \mathcal{D}^{*}$ (see Lemma 4(a) in Appendix \ref{Appendix:ALP}). According to Theorem \ref{Theorem:context_duality}(b), $\mathcal{L}H^{LP}(\textbf{x})\leq H^{LP}(\textbf{x})$. 
\item[(c)]  Suppose $H(\textbf{x}_{0})- (\textbf{R}(\textbf{x}_0,\textbf{a}_0)+\beta \mathbb{E}[ H(\textbf{x}_{1})|\textbf{x}_0,\textbf{a}_0] ) \geq  \varepsilon$ for all  $\textbf{x}_0\in \mathcal{X}$  and $\textbf{a}_0\in \bar{\mathcal{A}}(\textbf{x}_0)$. Then for all $\textbf{a}\in \mathbb{\bar{A}}$,
$$I_H(\textbf{a},\omega,\tau)=\sum_{t=0}^{\tau}\big(\textbf{R}( \textbf{x}_t,\textbf{a}_t)+\beta \mathbb{E}[ H(\textbf{x}_{t+1})|\textbf{x}_t,\textbf{a}_t]- H(\textbf{x}_{t})\big)\leq -\varepsilon(\tau+1),$$
for any $\tau\in \mathbb{N}$ and $\omega\in \Omega$. Therefore, $\mathbb{E}[\max_{\textbf{a}\in \mathcal{\bar{A}}(\tau)}\{I_H(\textbf{a},\omega,\tau)\}|\textbf{x}_0]\leq \mathbb{E}_0[-\varepsilon(\tau+1)]= \frac{-\varepsilon }{1-\beta},$ which implies $\mathcal{L}H(\textbf{x}) \leq H(x_{0})- \frac{\varepsilon}{1-\beta}.$
\end{enumerate}
\end{IEEEproof}

The last condition is used to  measure the gap between $H$ and $\mathcal{L}H(\textbf{x})$. A natural question is whether the improvement of the information relaxation bound over the Lagrangian bound can be significant. In Appendix \ref{Section:Information_Bound:Example}, we provide an affirmative answer by investigating the example proposed in \cite{adelman2008relaxations}, where the Lagrangian bound can be arbitrarily poor compared with the optimal value; as opposed to the performance of the Lagrangian bound, we show that  the optimal value can be recovered by improving the Lagrangian bound by the information relaxation approach.

A significant difference of the information relaxation and Lagrangian relaxation methods in the weakly coupled dynamic program is that the strong duality exists  in the former relaxation (at least theoretically), while such a result does not hold in general for the latter approach. The following theorem characterizes  the sufficient and necessary conditions such that $V(\textbf{x};\alpha')=\mathcal{L}{H}(\textbf{x}_0)$, where  $\alpha'\in \mathbb{\bar{A}}_{\mathbb{F}}$ is  a stationary Markov policy and $H$ is a function in $\mathcal{D}$.
\begin{theorem}\label{Theorem:Information_Complementary}
  Let $H\in \mathcal{D}$ and a stationary Markov policy  $\alpha'=(\alpha'_{\delta},\alpha'_{\delta},\cdots)\in \mathbb{\bar{A}}_{\mathbb{F}}$ , i.e., $\alpha'_{\delta}(\textbf{x})\in \mathcal{\bar{A}}(\textbf{x})$. A necessary and sufficient condition for  $V(\textbf{x}_0;\alpha')=\mathcal{L}H(\textbf{x}_0) $ for all $\textbf{x}_0\in \mathcal{X}$ is  that
\begin{align}
&\max_{\textbf{a}\in \mathcal{\bar{A}}(T)} \left\{\sum_{t=0}^{T}\big(R( \textbf{x}_t,\textbf{a}_t)+\beta \mathbb{E}\left[ H(\textbf{x}_{t+1})|\textbf{x}_t,\textbf{a}_t\right]- H(\textbf{x}_{t})\big)\right\} \notag\\
=&\sum_{t=0}^{T}\big(R(\textbf{x}_t,\alpha'_{\delta}(\textbf{x}_t))+\beta \mathbb{E}[ H(\textbf{x}_{t+1})|\textbf{x}_t,\alpha'_{\delta}(\textbf{x}_t)]- H(\textbf{x}_{t})\big) \label{context_complementary21}
\end{align}
for  $\omega\in \Omega$ almost surely,  $T=0,1,2,\cdots.$
In particular, by considering the case $T=0$,
$$\alpha'_{\delta}(\textbf{x}_0)\in \arg \max_{\textbf{a}_0\in  \mathcal{\bar{A}}(\textbf{x}_0)} \{R(\textbf{x}_0,\textbf{a}_0)+\beta \mathbb{E}[ H(\textbf{x}_{1})|\textbf{x}_0,\textbf{a}_0]\} $$
for all $\textbf{x}_0\in \mathcal{X}.$
\end{theorem}

The proof of Theorem\ref{Theorem:Information_Complementary} is in Appendix \ref{Appendix:ProofTheorem4}. Theorem\ref{Theorem:Information_Complementary} characterizes the optimality conditions of a  policy $\alpha'$ to (\ref{context_value_function}) and value approximation $H$ in (\ref{context_dual_rep}) as a pair: the optimal policy to the inner optimization problem of any horizon $T$ induced by the approximate value function   is non-anticipative and also stationary, though these decisions can be  chosen to be anticipative and non-stationary. In particular, the  policy $\alpha'$  is equal to a greedy policy  induced by the approximate value function $H$.

We connect Theorem \ref{Theorem:Information_Complementary} to  the analogous conditions for the Lagrangian bound  in  Theorem 2 and Lemma 1 of \cite{adelman2008relaxations}.  We review  the  sufficient and necessary conditions therein.  To ease comparison, we present them in a parallel way as the statement of  Theorem \ref{Theorem:Information_Complementary}.
\begin{lemma}\label{Propostion:Lagrangian_Complementary}
Let $\boldsymbol\lambda^{\circ}\geq 0$ and a stationary Markov policy  $\alpha^{\circ}=(\alpha^{\circ}_{\delta},\alpha^{\circ}_{\delta},\cdots)\in \mathbb{\bar{A}}_{\mathbb{F}}$, i.e., $\alpha^{\circ}_{\delta}(\textbf{x})\in \bar{\mathcal{A}}(\textbf{x})$. A necessary and sufficient condition for  $V(\textbf{x}_0;\alpha^{\circ})=J^{\boldsymbol\lambda^{\circ}}(\textbf{x}_0)$  for all $\textbf{x}_0\in \mathcal{X}$ is  that for all $\textbf{x}_0\in \mathcal{X}$, $\boldsymbol\lambda^{\circ \top}\left[\textbf{b}-\textbf{B}(\textbf{x}_0,\alpha^{\circ}_{\delta}(\textbf{x}_0))\right]=0$ and
\begin{align}
\alpha^\circ_{ \delta}(\textbf{x}_0)\in \arg\max_{\textbf{a}_0\in \mathcal{A}(\textbf{x}_0)}\left\{ \textbf{R}(\textbf{x}_0,\textbf{a}_0)+\boldsymbol\lambda^{\circ \top}[\textbf{b}-\textbf{B}(\textbf{x}_0,\textbf{a}_0)]+\beta
\mathbb{E}[J^{\boldsymbol\lambda^{\circ}}(\textbf{x}_1)|\textbf{x}_0,\textbf{a}_0]\right\}.
 \label{context_complementary11}
\end{align}
\end{lemma}

The conditions in Lemma  \ref{Propostion:Lagrangian_Complementary} are more stringent than those in Theorem \ref{Theorem:Information_Complementary}, since for any $\boldsymbol\lambda^{\circ}\geq 0$ and $\alpha^{\circ}=(\alpha^{\circ}_{\delta},\alpha^{\circ}_{\delta},\cdots)\in \mathbb{\bar{A}}_{\mathbb{F}}$,
$$V(\textbf{x}_0;\alpha^{\circ})\leq \mathcal{L}J^{\boldsymbol\lambda^{\circ}}(\textbf{x}_0)\leq\mathcal{L}^{\boldsymbol\lambda^{\circ}}J^{\boldsymbol\lambda^{\circ}}(\textbf{x}_0)= J^{\boldsymbol\lambda^{\circ}}(\textbf{x}_0)~~\text{for all}~\textbf{x}_0\in \mathcal{X}.$$

We show the connection of  Theorem \ref{Theorem:Information_Complementary} to Lemma \ref{Propostion:Lagrangian_Complementary} in the following.
If  $V(\textbf{x}_0;\alpha^{\circ})=J^{\boldsymbol\lambda^{\circ}}(\textbf{x}_0)$ for some $\alpha^{\circ}\in \mathbb{\bar{A}}_{\mathbb{F}}$ and $\boldsymbol\lambda^{\circ}\geq 0$, it implies $\mathcal{L}J^{\boldsymbol\lambda^{\circ}}(\textbf{x}_0)=\mathcal{L}^{\boldsymbol\lambda^{\circ}}J^{\boldsymbol\lambda^{\circ}}(\textbf{x}_0)$ and $V(\textbf{x}_0;\alpha^{\circ})=\mathcal{L}J^{\boldsymbol\lambda^{\circ}}(\textbf{x}_0)$.   Therefore,  the inequality (\ref{WeaklyCoupled_pathwise_inequality}) is actually an equality implied by $\mathcal{L}J^{\boldsymbol\lambda^{\circ}}(\textbf{x}_0)=\mathcal{L}^{\boldsymbol\lambda^{\circ}}J^{\boldsymbol\lambda^{\circ}}(\textbf{x}_0)$:
\begin{align*}
&\max_{\textbf{a}\in \mathcal{A}(\tau)} \left\{\sum_{t=0}^{T}\left(R(\textbf{x}_t,\textbf{a}_t)+\boldsymbol\lambda^{\circ\top}[\textbf{b}-\textbf{B}(\textbf{x}_t,\textbf{a}_t)]+\beta \mathbb{E}[ J^{\boldsymbol\lambda^{\circ}}(\textbf{x}_{t+1})|\textbf{x}_t,\textbf{a}_t]- J^{\boldsymbol\lambda^{\circ}}(\textbf{x}_t)\right)\right\}\\
=&\max_{\textbf{a}\in \mathcal{\bar{A}}(\tau)} \left\{\sum_{t=0}^{T}\big(R( \textbf{x}_t,\textbf{a}_t)+\beta \mathbb{E}\left[ J^{\boldsymbol\lambda^{\circ}}(\textbf{x}_{t+1})|\textbf{x}_t,\textbf{a}_t\right]- J^{\boldsymbol\lambda^{\circ}}(\textbf{x}_t)\big)\right\} \notag\\
=&\sum_{t=0}^{T}\big(R(\textbf{x}_t,\alpha'_{\delta}(\textbf{x}_t))+\beta \mathbb{E}[ J^{\boldsymbol\lambda^{\circ}}(\textbf{x}_{t+1})|\textbf{x}_t,\alpha'_{\delta}(\textbf{x}_t)]- J^{\boldsymbol\lambda^{\circ}}(\textbf{x}_t)\big),
\end{align*}
where the second equality holds for any $T\in \mathbb{N}$ due to $V(\textbf{x}_0;\alpha^{\circ})=\mathcal{L}J^{\boldsymbol\lambda^{\circ}}(\textbf{x}_0)$, implied by Theorem 3. Consider the special case $T=0$ and recall that $\boldsymbol\lambda^{\circ \top}\left[\textbf{b}-\textbf{B}(\textbf{x}_0,\alpha^{\circ}_{\delta}(\textbf{x}_0))\right]\geq0$, it is simple to verify  (\ref{context_complementary11}) and the condition $\boldsymbol\lambda^{\circ \top}\left[\textbf{b}-\textbf{B}(\textbf{x}_0,\alpha^{\circ}_{\delta}(\textbf{x}_0))\right]=0$ in Lemma \ref{Propostion:Lagrangian_Complementary}.



\section{Practical Information Relaxation Bound for Large-scale Problems }\label{Section:Computational_Method}
The information relaxation approach has the desirable property that it generates tighter upper bound based on the Lagrangian bound; however,  computing the  information relaxation bound can be   challenging  in large-scale weakly coupled dynamic program due to the  intractable  inner optimization problem. To be specific, the size of this scenario-dependent  optimization problem increases exponentially with respect to the number of the projects or subproblems $N$, and also increases at least linearly in the  horizon  $\tau$. Instead of computing the optimal value of the inner optimization problem, we discuss how to  derive its upper bound that is computationally tractable. Therefore, this sub-optimal method still leads to a valid upper bound on the value function, which is referred to as the ``practical information relaxation bound''. We will show its performance guarantee under certain conditions.


Throughout this section we assume that the approximate value function is of the additively separable form $H(\textbf{x})=\theta+\sum_{n=1}^NH^n(x^n)$, where $\theta$ is a constant and $H^n:\mathcal{X}^n\rightarrow \mathbb{R}$ for $n=1,\cdots,N$. We denote by $\mathcal{D}^\circ$ the space of additively separable functions. By substituting $H(\cdot)$ in   (\ref{context_inneropt}) by $\theta+\sum_{n=1}^NH^n(\cdot)$, we can rewrite $I_H$ as
\begin{equation}
I_{H}(\textbf{a},\omega,\tau)=\sum_{n=1}^N\left[\sum_{t=0}^{\tau} \left(R^n(x^n_t,a^n_t)+\beta \mathbb{E}[ H^n(x^n_{t+1})|x^n_t,a^n_t]- H^n(x^n_{t})\right)\right]-(\tau+1)(1-\beta)\theta . \label{I_H}
\end{equation}

\subsection{Relaxation of the Inner Optimization Problem }
Noting that  the scenario-dependent inner optimization problem $\max_{\textbf{a}\in \bar{\mathcal{A}}(\tau)}\{I_{H}(\textbf{a},\omega,\tau)\}$ is also \emph{weakly coupled} due to the additively separable structure of (\ref{I_H}) and the feasible control set $\bar{\mathcal{A}}(\tau)$.  To obtain an upper bound on its optimal value,   we  dualize the linking constraints for each period up to time $\tau$, and introduce the Lagrangian function $I_H(\textbf{a},\omega,\tau;\boldsymbol\mu)$ for $\boldsymbol\mu\triangleq(\boldsymbol\mu_0,\cdots,\boldsymbol\mu_{\tau})$ with each $\boldsymbol\mu_t\in \mathbb{R}^{L}_{+}$:
\begin{align}
I_{H}(\textbf{a},\omega,\tau;\boldsymbol\mu)\triangleq& \sum_{t=0}^{\tau}\left(\textbf{R}(\textbf{x}_t,\textbf{a}_t)+\beta \mathbb{E}[ H(\textbf{x}_{t+1})|\textbf{x}_t,\textbf{a}_t]- H(\textbf{x}_{t})+\boldsymbol\mu_t^{\top}[\textbf{b}-\textbf{B}(\textbf{x}_t,\textbf{a}_t)]\right) \notag\\
=& \sum_{t=0}^{\tau}\bigg[ \sum_{n=1}^N \big(R^n(x^n_t,a^n_t)+\beta \mathbb{E}[ H^n(x^n_{t+1})|x^n_t,a^n_t]- H^n(x^n_{t})-\boldsymbol\mu_t^{\top}\textbf{B}^n(x^n_t,a^n_t)\big) \notag\\
&-(1-\beta)\theta\bigg] + \sum_{t=0}^{\tau}\boldsymbol \mu_t^{\top}\textbf{ b} \notag\\
=&\sum_{n=1}^N I^n_{H^n}(\textbf{a}^n,\omega,\tau;\boldsymbol\mu) -(\tau+1)(1-\beta)\theta  + \sum_{t=0}^{\tau} \boldsymbol\mu_t^{\top} \textbf{b}, \label{separable}
\end{align}
where  $I^n_{H^n}$ in (\ref{separable}) is defined as
$$I^n_{H^n}(\textbf{a}^n,\omega,\tau;\boldsymbol\mu)\triangleq\sum_{t=0}^{\tau}\big( R^n(x^n_t,a^n_t)+\beta \mathbb{E}[ H^n(x^n_{t+1})|x^n_t,a^n_t]- H^n(x^n_{t})-\boldsymbol\mu_t^{\top}\textbf{B}^n(x^n_t,a^n_t)\big)$$
with $\textbf{a}^n\triangleq(a^n_0,\cdots,a^n_{\tau})$. In particular, $I_{H}(\textbf{a},\omega,\tau)=I_{H}(\textbf{a},\omega,\tau;0)$. Given any $\boldsymbol\mu\geq0$,  it is straightforward to see
$$\max_{\textbf{a}\in \bar{\mathcal{A}}(\tau)} I_{H}(\textbf{a},\omega,\tau)\leq \max_{\textbf{a}\in \mathcal{A}(\tau)}I_H(\textbf{a},\omega,\tau;\boldsymbol\mu).$$

According to (\ref{separable}), the Lagrangian dual $\max_{\textbf{a}\in \mathcal{A}(\tau)}\{I_H(\textbf{a},\omega,\tau;\boldsymbol\mu)\}$ can be decomposed as
\begin{align}
\max_{\textbf{a}\in\mathcal{A}(\tau)} \{I_{H}(\textbf{a},\omega,\tau;\boldsymbol\mu)\}=\sum_{n=1}^N \max_{\textbf{a}^n\in \mathcal{A}^n(\tau)}\left\{I^n_{H^n}(\textbf{a}^n,\omega,\tau;\boldsymbol\mu)\right\} -(\tau+1)(1-\beta)\theta  + \sum_{t=0}^{\tau}\boldsymbol \mu_t^{\top} \textbf{b}, \label{inner_opt_decomp}
\end{align}
where $\mathcal{A}^n(\tau)\triangleq\mathcal{A}^n(x^n_0)\times \cdots\times \mathcal{A}^n(x^n_{\tau})$. The equality (\ref{inner_opt_decomp}) implies that the computational cost on solving $\max_{\textbf{a}\in\mathcal{A}(\tau)} \{I_{H}(\textbf{a},\omega,\tau;\boldsymbol\mu)\}$  is linear rather than exponential in the number of the subproblems $N$. Therefore, the Lagrangian relaxation significantly reduces the computational complexity, and hence solving (\ref{inner_opt_decomp}) to optimality becomes potentially tractable.

It remains to find the optimal $\boldsymbol\mu^*$ that achieves the minimum of $I_{H}(\textbf{a},\omega,\tau;\boldsymbol\mu)$ over $\boldsymbol\mu\geq 0$. To this end, we list some properties of  $\max_{\textbf{a}\in \mathbb{A}(\tau)} I_{H}(\textbf{a},\omega,\tau;\boldsymbol\mu)$ as a function of $\boldsymbol\mu$,  based on properties of  Lagrangian relaxation.

%

\begin{lemma} \label{Lemma:Inner_Lagragian_relaxation} Given  $I_{H}(\textbf{a},\omega,\tau;\boldsymbol\mu)$  defined in (\ref{separable}), where $\omega\in \Omega$ and $\tau\in \mathbb{N}$. Then
\begin{enumerate}
\item[(a)] $\max_{\textbf{a}\in \mathcal{A}(\tau)}I_{H}(\textbf{a},\omega,\tau;\boldsymbol\mu)$ is  convex in $\boldsymbol\mu$.

\item[(b)]Let $\textbf{a}^\circ =(\textbf{a}^\circ_0,\cdots, \textbf{a}^\circ_{\tau})\in\arg\max_{\textbf{a}\in \mathcal{A}(\tau)}I_{H}(\textbf{a},\omega,\tau;\boldsymbol\mu)$ for a fixed $\boldsymbol\mu\geq0$.  Then 
$$\left[\textbf{b}-\textbf{B}(\textbf{x}^\circ_0,\textbf{a}^\circ_0),\cdots,\textbf{b}-\textbf{B}(\textbf{x}^\circ_\tau,\textbf{a}^\circ_\tau)\right] \in \partial I_{H}(\textbf{a}^\circ,\omega,\tau;\boldsymbol\mu),$$
where $\textbf{x}^\circ_t=\textbf{x}_t(\textbf{a}^\circ(t-1),\omega)$ is the state trajectory under $\textbf{a}^\circ$ and $\omega$, and  $\partial I_{H}(\textbf{a}^\circ,\omega,\tau;\boldsymbol\mu)$ is the subdifferential of $I_{H}(\textbf{a},\omega,\tau;\boldsymbol\mu)$ with respect to $\boldsymbol\mu$ at $\textbf{a}=\textbf{a}^\circ$.

\item[(c)] $\max_{\textbf{a}\in \bar{\mathcal{A}}(\tau)} I_{H}(\textbf{a},\omega,\tau)\leq \min_{\boldsymbol\mu\geq 0} \max_{\textbf{a}\in \mathcal{A}(\tau)}I_{H}(\textbf{a},\omega,\tau;\boldsymbol\mu).$
\end{enumerate}
\end{lemma}

Lemma \ref{Lemma:Inner_Lagragian_relaxation} indicates that $\min_{\boldsymbol\mu\geq 0} \max_{\textbf{a}\in \mathcal{A}(\tau)}I_{H}(\textbf{a},\omega,\tau;\boldsymbol\mu)$  is a convex optimization problem in $\boldsymbol\mu$. Since its subgradient at $\boldsymbol\mu$ is known, we can employ the standard subgradient method or its variant to locate the optimal solution efficiently. Due to  Lemma \ref{Lemma:Inner_Lagragian_relaxation}(c), we refer to $\min_{\boldsymbol\mu\geq 0} \max_{\textbf{a}\in \mathcal{A}(\tau)}I_{H}(\textbf{a},\omega,\tau;\boldsymbol\mu)$ as the ``relaxed inner optimization problem''.

Based on the relaxed inner optimization problem we  define  a new operator $\mathcal{L}^\circ$ that can be viewed as a ``relaxed'' version of  $\mathcal{L}$ on the additively separable function space $\mathcal{D}^\circ$:
\begin{align}
\mathcal{L}^\circ H(\textbf{x}) \triangleq H(\textbf{x})+ \mathbb{E}_0\left[\min_{\boldsymbol\mu\geq 0} \max_{\textbf{a}\in \mathcal{A}(\tau)}I_{H}(\textbf{a},\omega,\tau;\boldsymbol\mu)\right]. \label{context_dual_rep_practical}
\end{align}
Due to the computational tractability of $\mathcal{L}^\circ H(\textbf{x})$, it will be referred to as  ``practical information relaxation bound''. In the next theorem we formalize the bound performance of $\mathcal{L}^\circ H(\textbf{x})$, which naturally places an upper bound on the  information relaxation bound $\mathcal{L}H$; moreover, the performance of $\mathcal{L}^\circ  J^{\boldsymbol\lambda}(\textbf{x})$ is no worse than the Lagrangian bound $J^{\boldsymbol\lambda}(\textbf{x})$.
\begin{theorem}\label{Theorem:Relaxed_inner_opt}
Suppose $H\in \mathcal{D}^\circ$.  Then
\begin{enumerate}
\item[(a)] $\mathcal{L}H(\textbf{x})\leq \mathcal{L}^\circ H(\textbf{x})$ for all $\textbf{x}\in \mathcal{X}$.
\item[(b)] Suppose  $H=J^{\boldsymbol\lambda}$ is a Lagrangian bound for some $\boldsymbol\lambda\geq0$. Then for every $\omega\in \Omega$ and $\tau\in \mathbb{N}$,
$$\min_{\boldsymbol\mu\geq 0} \max_{\textbf{a}\in \mathcal{A}(\tau)}\left\{I_{J^{\boldsymbol\lambda}}(\textbf{a},\omega,\tau;\boldsymbol\mu)\right\}\leq 0.$$
 Consequently, $\mathcal{L}^\circ J^{\boldsymbol\lambda}(\textbf{x})\leq  J^{\boldsymbol\lambda}(\textbf{x})$ for all $\textbf{x}\in \mathcal{X}$.
 \end{enumerate}
\end{theorem}
\begin{IEEEproof}
\begin{enumerate}
\item[(a)] This is because for every  $\omega\in \Omega$ and $\tau\in \mathbb{N}$,
$$\max_{\textbf{a}\in \bar{\mathcal{A}}(\tau)} \{I_{H}(\textbf{a},\omega,\tau)\}\leq \min_{\boldsymbol\mu\geq 0} \max_{\textbf{a}\in \mathcal{A}(\tau)}\{I_{H}(\textbf{a},\omega,\tau;\boldsymbol\mu)\}.$$
\item[(b)] Note that $J^{\boldsymbol\lambda}(\textbf{x}_0)=\mathcal{L}^{\boldsymbol\lambda}J^{\boldsymbol\lambda}(\textbf{x}_0)$ and $\mathcal{L}^{\boldsymbol\lambda}J^{\boldsymbol\lambda}(\textbf{x}_0)=J^{\boldsymbol\lambda}(\textbf{x}_0)+\mathbb{E}_0\left[\max_{\textbf{a}\in\mathcal{A}(\tau)} \{I_{J^{\boldsymbol\lambda}}(\textbf{a},\omega,\tau;\boldsymbol\lambda)\}\right]$ (see the definition of $\mathcal{L}^{\boldsymbol\lambda}J^{\boldsymbol\lambda}$ in (\ref{Lagragian_duality})). Given any $\boldsymbol\lambda\geq 0$,  we have for every $\omega\in \Omega$ and $\tau\in \mathbb{N}$,
\begin{align}
0=\max_{\textbf{a}\in \mathcal{A}(\tau)} \{  I_{J^{\boldsymbol\lambda}}(\textbf{a},\omega,\tau;\boldsymbol\lambda)  \}\geq \min_{\boldsymbol\mu\geq 0} \max_{\textbf{a}\in \mathcal{A}(\tau)} \{I_{J^{\boldsymbol\lambda}}(\textbf{a},\omega,\tau;\boldsymbol\mu)\},   \label{inner_compare}
\end{align}
where the first equality follows (\ref{inner_LagrangianRelaxation}) in Lemma \ref{Corollary:Dual_Lagrangian}. Therefore,
\begin{align*}
\mathcal{L}^{\boldsymbol\lambda}J^{\boldsymbol\lambda}(\textbf{x}_0)\geq J^{\boldsymbol\lambda}(\textbf{x}_0)+\mathbb{E}_0\left[\min_{\boldsymbol\mu\geq 0}\max_{\textbf{a}\in \mathcal{A}(\tau)} \left\{I_{J^{\boldsymbol\lambda}}(\textbf{a},\omega,\tau;\mu)\right\}\right]=\mathcal{L}^\circ J^{\boldsymbol\lambda}(\textbf{x}_0).
\end{align*}
\end{enumerate}
\end{IEEEproof}

The inequality (\ref{inner_compare}) highlights the  comparison between two scenario-based inner optimization problems: the right  term of the inequality in (\ref{inner_compare}) allows $\boldsymbol\mu=(\boldsymbol\mu_0,\cdots,\boldsymbol\mu_{\tau})$ \big(contained in $\sum_{t=0}^{\tau}\boldsymbol\mu_t^{\top}[\textbf{b}-\textbf{B}(\textbf{x}_t,\textbf{a}_t)]$\big)   to be different across  periods; on the other hand, the left term forces
$\boldsymbol\mu=(\boldsymbol\lambda,\cdots,\boldsymbol\lambda)$ \big(contained in $\sum_{t=0}^{\tau}\boldsymbol\lambda^{\top}[\textbf{b}-\textbf{B}(\textbf{x}_t,\textbf{a}_t)]$\big) to be constant over time. Therefore, $\mathcal{L}^\circ J^{\boldsymbol\lambda}$ can be viewed as an intermediate relaxation between the ``exact'' information relaxation $\mathcal{L}J^{\boldsymbol\lambda}$ and the Lagrangian relaxation $J^{\boldsymbol\lambda}(=\mathcal{L}^{\boldsymbol\lambda}J^{\boldsymbol\lambda})$. Another useful observation is that $\boldsymbol\mu=(\boldsymbol\lambda,\cdots,\boldsymbol\lambda)$  can naturally serve as the initial point to solve $\min_{\boldsymbol\mu\geq 0} \max_{\textbf{a}\in \mathcal{A}(\tau)}\{I_{J^{\boldsymbol\lambda}}(\textbf{a},\omega,\tau;\boldsymbol\mu)\}$ via the subgradient method.


Note that the computational complexity of the inner optimization problem also depends on the time horizon $\tau$. In case of drawing a sample of $\tau$ that is a large number (often occurs when $\beta$ that is close to $1$), we propose a simple remedy to ease  computation, i.e., to  truncate the random horizon of the relaxed inner optimization problem up to some deterministic time $\mathcal{T}\in \mathbb{N}$ that is sufficiently large.  This operation reduces the computational cost in some extreme cases; however, we would also like to know the quality of the resulting bounds with regard to different $\mathcal{T}'$s.  The following result shows  the complexity versus quality trade-off in choosing an appropriate $\mathcal{T}$: a greater truncated horizon $\mathcal{T}$ implies a more difficult inner optimization problem but  guarantees  better bound.

\begin{corollary} \label{Corollary:Truncation}
Suppose $\mathcal{T}\in \mathbb{N}$. Define
$$\mathcal{L}^{\circ}_{\mathcal{T}}J^{\boldsymbol\lambda}(\textbf{x})\triangleq J^{\boldsymbol\lambda}(\textbf{x})+ \mathbb{E}_0\left[\min_{\boldsymbol\mu\geq 0} \max_{\textbf{a}\in \mathcal{A}(\tau)}I_{J^{\boldsymbol\lambda}}(\textbf{a},\omega,\tau\wedge \mathcal{T};\boldsymbol\mu)\right],$$
where $\tau \wedge \mathcal{T} =\min\{\tau,\mathcal{T}\}$. Then
\begin{enumerate}
\item[(a)] $\mathcal{L}^\circ J^{\boldsymbol\lambda}(\textbf{x})\leq \mathcal{L}^{\circ}_{\mathcal{T}+1}J^{\boldsymbol\lambda}(\textbf{x}) \leq \mathcal{L}^{\circ}_{\mathcal{T}}J^{\boldsymbol\lambda}(\textbf{x})\leq J^{\boldsymbol\lambda}(\textbf{x}).$
\item[(b)] $\lim_{\mathcal{T}\rightarrow\infty} \mathcal{L}^{\circ}_{\mathcal{T}}J^{\boldsymbol\lambda}(\textbf{x})=\mathcal{L}^\circ J^{\boldsymbol\lambda}(\textbf{x}).$
\end{enumerate}
\end{corollary}
\begin{IEEEproof}
Note that by fixing $\omega\in\Omega$ and $\tau\in \mathbb{N}$, the following inequality holds for any  $\mathcal{T}\in \mathbb{N}$:
$$\min_{\boldsymbol\mu\geq 0} \max_{\textbf{a}\in \mathcal{A}(\tau)}I_{J^{\boldsymbol\lambda}}(\textbf{a},\omega,\tau;\boldsymbol\mu)\leq \min_{\boldsymbol\mu\geq 0} \max_{\textbf{a}\in \mathcal{A}(\tau)}I_{J^{\boldsymbol\lambda}}(\textbf{a},\omega,\tau\wedge (\mathcal{T}+1);\boldsymbol\mu)\leq  \min_{\boldsymbol\mu\geq 0} \max_{\textbf{a}\in \mathcal{A}(\tau)}I_{J^{\boldsymbol\lambda}}(\textbf{a},\omega,\tau\wedge \mathcal{T};\boldsymbol\mu)\leq 0.$$
 Therefore, the inequality in $(a)$ follows from the above inequality immediately, and the equality in $(b)$ is true due to the monotone convergence theorem.
\end{IEEEproof}

\subsection{The Gap between Practical and Exact Information Relaxation Bounds}
The practical information relaxation bound $\mathcal{L}^\circ H(x)$ effectively reduces the computational cost compared to deriving the exact information relaxation bound $\mathcal{L}H(x)$, though yields a less tight bound.  In this subsection we investigate  the gap  $\mathcal{L}^\circ H(x) -\mathcal{L} H(x)$, which is the average difference between the optimal values of the exact and relaxed  inner optimization problems, i.e.,
\begin{equation}
\min_{\boldsymbol\mu\geq 0} \max_{\textbf{a}\in \mathcal{A}(\tau)}I_{H}(\textbf{a},\omega,\tau;\boldsymbol\mu)-\max_{\textbf{a}\in \mathcal{\bar{A}}(\tau)} I_{H}(\textbf{a},\omega,\tau).\label{gap_inner_opt}
\end{equation}
\cite{bertsekas1982constrained} established the sufficient conditions such that the Lagrangian duality gap of  the weakly coupled deterministic optimization problem is uniformly bounded regardless of the number of the subproblems (see Appendix \ref{Appendix:Separable_DulaityGap}). We will show a similar result for $\mathcal{L}^\circ H(x) -\mathcal{L} H(x)$  assuming that $H$ is additively separable.

 We begin with an intuitive interpretation  on the duality gap (\ref{gap_inner_opt})  by looking at two equivalent linear program formulations of  (\ref{inner_opt_decomp}). We fix  $\omega\in \Omega$ and $\tau\in \mathbb{N}$, and assume that  the control space $\mathcal{A}$ is finite. For each project $n=1,\cdots,N$, we can then  enumerate all  state trajectories  of $(x^n_1,\cdots,x^n_\tau)$  (denoted by $(x^{n,n_k}_1,\cdots,x^{n,n_k}_\tau)$ with index $n_k$) associated with the control sequence $(a^n_1,\cdots,a^n_\tau)\in \mathcal{A}^n(\tau)$ (denoted by $\textbf{a}^{n,n_k}$). Noting that  $I^n_{H}(\textbf{a}^{n,n_k},\omega,\tau;0)-\sum_{t=0}^{\tau} \boldsymbol\mu_t ^{\top}\textbf{ B}^n_t(x^{n,n_k}, a^{n,n_k})= I^n_{H}(\textbf{a}^{n,n_k},\omega,\tau;\boldsymbol\mu) $. Then  (\ref{inner_opt_decomp}) can be equivalently written as the following linear program,
\begin{align}
\min_{ \{y_n, \boldsymbol\mu_t\}} &\sum_{n=1}^N y_n  + \sum_{t=0}^{\tau} \boldsymbol\mu_t^{\top} \textbf{b}- (\tau+1) (1-\beta)\theta \notag\\
\text{s.t.}~~ & y_n \geq I^n_{H}(\textbf{a}^{n,n_k},\omega,\tau;0)-\sum_{t=0}^{\tau} \boldsymbol\mu_t^{\top}\textbf{B}^n(x^{n,n_k}, a^{n,n_k})  ~~\text{for all}~ n_k,~  n=1,\cdots,N;  \label{LP:prob}\\
     &\boldsymbol\mu_t\geq 0,  ~~t=0,\cdots, \tau.  \notag 
\end{align}
We use $p^{n,n_k}$ to denote the dual variable associated with (\ref{LP:prob}), so the dual linear program is
\begin{align*}
\max_{ \{p^{n,n_k}\}} &\sum_{n=1}^N \sum_{n_k}p^{n,n_k}I^n_{H}(\textbf{a}^{n,n_k},\omega,\tau;0) -(\tau+1) (1-\beta)\theta \\
\text{s.t.}~~ &\sum_{n=1}^N\sum_{n_k}p^{n,n_k} \textbf{B}^n(x^{n,n_k}_t,a^{n,n_k}_t)\leq \textbf{b}, ~~t=1,\cdots, \tau;  \\
     &\sum_{n_k}p^{n,n_k}=1,  ~~n=1,\cdots, N;  \\
     &p^{n,n_k}\geq0 ~\text{for all}~n_k ~and ~n=1,\cdots,N,
\end{align*}
where $p^{n,n_k}$ can be interpreted as the probability assigned to the  $n_k$-th scenario associated with project $n$. Comparing  the above linear program to (\ref{I_H}), it can be seen that  the feasible control set  $\bar{\mathcal{A}}(\tau)$ is enlarged to include all the randomized controls subject to the linking constraint. Therefore, the relaxed inner optimization problem can be viewed as the convexification of the exact inner optimization problem.  In addition, the optimal solution to  the dual linear program (that can be found by LP solvers at least for  finite  $|\mathcal{A}|$ and small horizon $\tau$) also provides benchmark result on (\ref{inner_opt_decomp}), which may help to adjust the parameters used in the subgradient method.

To characterize the gap  $\mathcal{L}^\circ H(x) -\mathcal{L} H(x),$  we list some technical assumptions based on Lemma 2 in Appendix \ref{Appendix:Separable_DulaityGap}. In particular, we denote $\textbf{B}^n(x^n_t,a^n_t)$ equivalently as $\textbf{B}^n_t(\textbf{a}^n,\omega)$, as $x_t^n$ depends on $\textbf{a}^n$ and $\omega$.
\begin{assumption}
For every state $\textbf{x}\in\mathcal{X}$, $\bar{\mathcal{A}}(\textbf{x})\neq \phi$.
 \end{assumption}

\begin{assumption}
Given $\omega\in \Omega$ and $T\in \mathbb{N}$, the sets
 $$S_n\triangleq\left\{(\textbf{a}^n, \textbf{B}^n_0(\textbf{a}^n,\omega),\cdots, \textbf{B}^n_{T}(\textbf{a}^n,\omega), I^n_{H}(\textbf{a}^n,\omega,T) )|\textbf{a}^n\in \mathcal{A}^n(T)\right\}$$
 are non-empty and compact for $n=1,\cdots,N$.
\end{assumption}
This assumption is automatically true if each $\mathcal{A}^n$ is finite, or $\mathcal{A}^n(T)$ is compact and each $\textbf{B}^n_t(\textbf{a}^n,\omega)$ and $I^n_{H}(\textbf{a}^n,\omega,T)$ are continuous functions on $\mathcal{A}^n(T)$.

\begin{assumption}
Given $\omega\in \Omega$ and $T\in \mathbb{N}$. For every $n=1,\cdots,N$,  we assume that for any $\tilde{\textbf{a}}^n\in conv(\mathcal{A}^n(T))$, there exists $\textbf{a}^n\in \mathcal{A}^n(T)$ such that
\begin{equation}
\textbf{B}^n_t(\textbf{a}^n, \omega)\leq  (\check{\text{cl}}~\textbf{B}^n_t)(\tilde{\textbf{a}}^n, \omega), ~t=0,\cdots,T,  \label{Assump3_inequ}
\end{equation}
where $\check{\text{cl}}~\textbf{B}^n_t$ is the function whose component is the convex closure of the corresponding component of $\textbf{B}^n_t$, i.e.,
{\small
$$\check{\text{cl}}~\textbf{B}^n_t(\tilde{\textbf{a}}^n, \omega) \triangleq  \inf\left\{ \sum_{n_k}p^{n,n_k} \textbf{B}^n_t(\textbf{a}^{n,n_k}, \omega)\bigg| \tilde{\textbf{a}}^n= \sum_{n_k}p^{n,n_k} \textbf{a}^{n,n_k}, ~\textbf{a}^{n,n_k}\in \mathcal{A}^n(T); \sum_{n_k}p^{n,n_k}=1,~p^{n,n_k}\geq 0 \right\}.$$}
\begin{remark}
All the sums in the definition of $\check{\text{cl}}~\textbf{B}^n_t(\tilde{\textbf{a}}^n, \cdot)$  are finite sums.
\end{remark}
\end{assumption}
This assumption is not trivially satisfied, as (\ref{Assump3_inequ}) can be a vector inequality. However, there are several cases that we can directly verify  Assumption 3 is true.
\begin{enumerate}
\item[Case 1.] Each $|\mathcal{A}^n|$  is finite,  the number of the linking constraints $L=1$ (therefore, each inequality in (\ref{Assump3_inequ}) is a scalar inequality), and each $\textbf{B}^n_t(\textbf{a}^n, \omega)$ (i.e., $\textbf{B}^n(x^n_t,a^n_t)$) only depends on $a^n_t$. A typical example is the restless bandit problem, in which the linking constraint is $\sum_{n=1}^N\textbf{B}^n(x^n_t,a^n_t)=\sum_{n=1}^N a^n_t=1$ with $a_t^n\in\{0,1\}$.
\item[Case 2.] If $\mathcal{A}^n(T)$ is convex, and the components of each $\textbf{B}^n_t(\textbf{a}^{n}, \omega)$  are convex over  $\mathcal{A}^n(T)$ for $t=0,\cdots,T$. Then $conv(\mathcal{A}^n(T))=\mathcal{A}^n(T)$, and $(\check{\text{cl}}~\textbf{B}^n_t)(\tilde{\textbf{a}}^n, \omega)=\textbf{B}^n_t(\tilde{\textbf{a}}^n, \omega)$.
\end{enumerate}

%

We  present our main result on the  gap $\mathcal{L}^\circ H(x) -\mathcal{L} H(x)$.


\begin{theorem} \label{Theorem:Bound_Gap}
Suppose that $H$ is  of the additively separable form  $H(\textbf{x})=\theta+\sum_{n=1}^N H^n(x^n)$, and Assumptions 1-3 hold for every $\omega\in \Omega$ and $T\in \mathbb{N}$. Then for all $\textbf{x}\in \mathcal{X}$,
  \begin{equation}
  \mathcal{L}^\circ H(\textbf{x})-\mathcal{L}H(\textbf{x})\leq \frac{(L-1)\beta+L+1}{(1-\beta)^2}  \max_{n=1,\cdots, N} \Gamma^n, \label{Gap_Bound}
  \end{equation}
where
\begin{align*}
\Gamma^n=&\sup_{x^n_0\in \mathcal{X}^n, a^n_0 \in \mathcal{A}^n(x^n_0)}\{ R^n(x^n_0,a^n_0)+\beta \mathbb{E}[ H^n(x^n_{1})|x^n_0,a^n_0]- H^n(x^n_{0}) \}\\
&-\inf_{x^n_0\in \mathcal{X}^n, a^n_0 \in \mathcal{A}^n(x^n_0)}\{ R^n(x^n_0,a^n_0)+\beta \mathbb{E}[ H^n(x^n_{1})|x^n_0,a^n_0]- H^n(x^n_{0})\}.
\end{align*}
\end{theorem}

The proof of Theorem \ref{Theorem:Bound_Gap} is in Appendix E. Theorem \ref{Theorem:Bound_Gap} not only characterizes the gap between $\mathcal{L}^\circ H(x)$ and $\mathcal{L}H(x)$, but also allows  controlling this gap  by restricting the feasible region of  $\{H^{n}(\cdot)\}_{n=1}^{N}$. To be specific, we can add to the linear program (\ref{Lagragian_LP}) or (\ref{ALP}) the following constraints on the  \emph{Bellman error} of each subproblem (i.e., $R^n(x^n_0,a^n_0)+\beta \mathbb{E}[ H^n(x^n_{1})|x^n_0,a^n_0]- H^n(x^n_{0})$) :
$$\Gamma^{n,2}~\geq R^n(x^n_0,a^n_0)+\beta \mathbb{E}[ H^n(x^n_{1})|x^n_0,a^n_0]- H^n(x^n_{0})\geq -\Gamma^{n,1},~\text{for all}~(x^n_0,a^n_0) ~\text{with}~ a^n_0\in \mathcal{A}^n(x^n_0), $$
where  $\Gamma^{n,1}$ and $\Gamma^{n,2}$ are two positive numbers for $n=1,\cdots,N$. Suppose that there is a feasible solution to  the linear program (\ref{Lagragian_LP}) or (\ref{ALP}), then $\mathcal{L}^\circ H(\textbf{x})-\mathcal{L}H(\textbf{x})$ can be bounded by   $\frac{(L-1)\beta+L+1}{(1-\beta)^2}\max_{n=1,\cdots,N}\{\Gamma^{n,1}+\Gamma^{n,2}\}$.  Note that  the greater   $\Gamma^{n,1}$ and $\Gamma^{n,2}$ are, the larger the feasible region of $\{H^{n}(\cdot)\}_{n=1}^N$ is, which implies a tighter bound $J^{\boldsymbol\lambda}(\textbf{x})$ or $H^{LP}(\textbf{x})$; they may be used  to generate tighter bounds $\mathcal{L}J^{\boldsymbol\lambda}(\textbf{x})$ or $\mathcal{L}H(\textbf{x})$ according to Theorem \ref{Propostion:Lagrangian_Information_Bound}. As a trade-off, the gap  between the practical information relaxation bound $\mathcal{L}^\circ H(\textbf{x})$ and the exact $\mathcal{L}H(\textbf{x})$ may be enlarged.

As a corollary, Theorem \ref{Theorem:Bound_Gap} indicates that the gap $\mathcal{L}^\circ H(x)-\mathcal{L}H(x)$ has a uniform bound in $N$, if the Bellman errors of individual subproblems (and hence  $\Gamma^n$) are uniformly bounded for all state-action pairs $\{(x^n_{0},a^n_0)\}$. Therefore, the relative gap   $\frac{\mathcal{L}^\circ H(\textbf{x})-\mathcal{L}H(\textbf{x})}{N} $ vanishes as $N$ goes to infinity.  We provide an instance in which $\{\Gamma^n\}_{n=1}^N$ are uniformly bounded with mild conditions on rewards and linking constraints.  

\begin{corollary}\label{Corollary:Example}
\item[(a)] If $\{\Gamma^n\}_{n=1}^N$ are uniformly bounded for all subproblems, then $\mathcal{L}^\circ H(\textbf{x})-\mathcal{L}H(\textbf{x})$ is also uniformly bounded with respect to the number of subproblems  $N$.
\item[(b)]  Let $H(\textbf{x})=J^{ \boldsymbol\lambda}(\textbf{x})= \frac{\boldsymbol\lambda^{\top}\textbf{b}}{1-\beta}  + \sum_{n=1}^{N}H^{n, \boldsymbol\lambda}(x^n)$ for some $ \boldsymbol\lambda\geq 0$. Suppose  there exists a constant $C>0$ such that  $\{|R^n|, |R^n-\boldsymbol\lambda^{\top}\textbf{B}^n|\}_{n=1}^N$ are uniformly bounded by $C$. Then
 $\{\Gamma^n\}_{n=1}^N$ are uniformly bounded by $\frac{4C}{1-\beta}$.
\end{corollary}
\begin{IEEEproof}
\begin{enumerate}
\item[(a)] The result directly follows from Theorem \ref{Theorem:Bound_Gap}.
\item[(b)]  It can be seen from (\ref{Subproblem_Bellman}) that $\{H^{n,\lambda}\}_{n=1}^N$ are uniformly bounded by $\frac{C}{1-\beta}$, since $|R^n-\boldsymbol\lambda^{\top}\textbf{B}^n|\leq C$. Therefore, for all $n=1,\cdots, N$,
$$\frac{2C}{1-\beta}\geq R^n(x^n_0,a^n_0)+\beta \mathbb{E}[H^{n,\lambda}(x^n_{1})|x^n_0,a^n_0]- H^{n,\lambda}(x^n_{0})\geq -\frac{2C}{1-\beta}~\text{for all}~(x^n_0,a^n_0) ~\text{with}~ a^n_0\in \mathcal{A}^n(x^n_0),$$
i.e., $\{\Gamma_n\}_{n=1}^N$ are uniformly bounded by $\frac{4C}{1-\beta}$.
\end{enumerate}
\end{IEEEproof}

In other words, if the optimal value is proportional to the number of the subproblems, i.e., $N C_1\leq V\leq N C_2$ for some $C_1,C_2>0$ (e.g.,  $C_1(1-\beta)\leq|R^n|\leq C_2(1-\beta)$ for all $n=1\cdots,N$), then the relative gap  $\frac{\mathcal{L}^\circ H(\textbf{x})-V(\textbf{x})}{V(\textbf{x})}$ converges to the relative gap $\frac{\mathcal{L}H(\textbf{x})-V(\textbf{x})}{V(\textbf{x})}$ as the number of subproblems $N$ increases.


\begin{remark}
All the results presented in Section \ref{Section:Computational_Method} and Section \ref{Section:Computation} have counterparts in the finite-horizon setting; we refer the readers to Appendix \ref{Appendix:FinteHorizon} for details.
\end{remark}

\section{Numerical Examples}\label{Section:Computation}
To investigate the empirical  performance of the information relaxation bound, we test our method on both discrete-state  and continuous-state weakly-coupled stochastic optimization problems: one is the standard restless bandit problem, and the other one is a linear quadratic control problem with a non-convex linking constraint.  We compare some heuristic policies with both the Lagrangian bound and the practical information relaxation bound. We gain some insight on the quality of these relaxation in terms of the number of the  subproblems.

\subsection{Restless Bandit}
We consider a standard restless bandit problem that consists of  $N$ projects, where each project $n$ can take one of a finite number $|\mathcal{X}^n|$ of states. For  each project there are two actions to choose from:  \emph{active} ($a^n_t=1$) and \emph{passive} ($a^n_t=0$).  At each time period $t$, exactly $1$ project is chosen (corresponding to the linking constraint $\sum_{n=1}^N a^n_t=1$), and  its state transits according to the \emph{active transition probability} and receive an \emph{active reward}; for the rest projects that are not chosen, their states also transit but according to the respective \emph{passive transition probability} and earns a respective \emph{passive reward} that is zero. Projects are selected sequentially  in order to maximize a discounted infinite-horizon reward.

{\footnotesize
\begin{table}[htpb]
  \centering \caption{Restless Bandit} \label{Table:Result1}
  \begin{tabular}{|c|c|cc|cc|c|cc|cc|} \hline
  \multicolumn{1}{|c} {} &\multicolumn{1}{|c} {} &\multicolumn{2}{|c} {Lag. Policy} &\multicolumn{2}{|c|} {PD Policy}  &\multicolumn{1}{c|} {Lag. Bound } &\multicolumn{2}{c|} {Info. Relaxation } &\multicolumn{2}{c|} { Duality Gap }\\ \hline
    { N } &$\beta$  &Value &S.E. &Value &S.E.    &Value   &Value  &S.E. &1  &2

    \rule[1ex]{0pt}{1em}\\\hline\hline
    $10$  &$0.90$  &$9.0241$     &$0.0018$    &$9.0797$    &0.0017      &$9.2971$     &9.1785      &0.0028       &1.09\% &54.6\%  \\
    $20$  &$0.90$  &$9.4037$     &$0.0002$    &$9.5779$    &0.0003      &$9.6747$     &9.6196      &0.0004       &0.44\% &56.9\%  \\
    $50$  &$0.90$  &$9.6244$     &$0.0000$    &$9.7346$    &0.0003      &$9.8336$     &9.7511      &0.0006       &0.17\% &83.3\%  \\
    $10$  &$0.95$  &$17.8418$    &$0.0038$    &$18.1502$   &0.0030      &$18.6041$    &18.3907     &0.0049       &1.33\% &47.0\%  \\
    $20$  &$0.95$  &$18.7189$    &$0.0005$    &$19.1487$   &0.0009      &$19.3272$    &19.2164     &0.0006      &0.35\% &62.1\%  \\
    $50$  &$0.95$  &$19.1705$    &$0.0002$    &$19.4720$   &0.0004      &$19.6740$    &19.5562     &0.0028       &0.43\% &58.3\%  \\
    $10$  &$0.98$  &$42.1093$    &0.0107      &$45.1822$   &$0.0098$    &$46.5247$    &45.9544     &0.0253       &1.71\% &42.8\%  \\
    $20$  &$0.98$  &$46.3824$    &0.0018      &$47.7935$   &$0.0017$    &$48.2850$    &47.9801     &0.0062      &0.39\% &62.0\%  \\
    $50$  &$0.98$  &$47.6881$    &0.0011     &$48.5806$   &$0.0019$    &$49.5797$    &49.2176     &0.0014       &1.31\% &36.2\%  \\
\hline
\end{tabular}
\end{table}}For each subproblem $n=1,\cdots,N$, we set the cardinality of its state space $|\mathcal{X}^n|=10$ in our numerical experiments. In Table \ref{Table:Result1} we list the numerical results and the corresponding parameters including the number of projects $N$, the discount factor $\beta$. We generate random instances of  active and passive transition probabilities for each subproblem. Active rewards are sampled from the uniform distribution on $[0,1]$. We compute the upper and lower bounds on $V(\textbf{x}_0)$, where the initial condition  $\textbf{x}_0=(1,1,\cdots,1)^{\top}$.  We first solve the Lagrangian bound --``Lag. Bound'' (i.e., $J^{\boldsymbol\lambda^*}$)  via the linear program (\ref{Lagragian_LP}), where we set  the  distribution $\upsilon(\cdot)$ to be uniform  over all states. Based on the ``Lag. Bound'' we compute the practical information relaxation bound -- ``Info. Relaxation'' via (\ref{context_dual_rep_practical}): we generate  $100$ scenarios ($\tau$ and $\omega$) and solve the associated relaxed inner optimization problems. Here we truncate the random time horizon $\tau$ up to  $\mathcal{T}=50$, $100$, and $150$ for $\beta=0.9$, $0.95$, $0.98$, respectively, i.e, we actually compute $\mathcal{L}^{\circ}_{\mathcal{T}}J^{\boldsymbol\lambda^*}(\textbf{x})$ and apply the subgradient method with at most $200$, $400$, and $1000$ iterations, respectively, or until the norm of the subgradient is exactly zero. The actual number of iterations mainly depends on the realization of $\tau$: the greater $\tau$ is, generally the more iterations are needed to attain convergence in the subgradient method.  To obtain two lower bounds --``Lag. Policy'' and ``PD policy'', we also generate $100$ scenarios, and apply the one-step greedy policy induced by the ``Lag. Bound'' and the ``primal-dual'' policy developed in \cite{bertsimas2000restless}, respectively. To demonstrate the performance of ``Info. Relaxation'', we report the relative duality gaps in two categories:
\begin{align*}
\text{Duality Gap 1}&=\frac{\text{``Info.~Relaxation''~--~``PD~Policy''}}{\text{``PD~Policy''}},  \\
\text{Duality Gap 2}&=\frac{\text{``Lag.~Bound''~--~``Info.~Relaxation''}}{\text{``Lag.~Bound''~--~``PD~Policy''}}.
\end{align*}
The ``Duality Gap 1'' shows that the  relative gaps between the best lower and upper bounds, and the ``Duality Gap 2'' reports the percentage of the reduced duality gap by comparing ``Info. Relaxation'' and ``Lag. Bound'' to ``PD~Policy''.


We observe that in our numerical experiments  the gap between  ``Lag. Policy'' and ``Lag. Bound'' generally increases as $\beta$ approaches 1, and they are relatively larger in the problem with smaller number of projects.  The lower bounds can be significantly  improved by applying the primal-dual policy in all cases. On the other hand,  ``Info. Relaxation'' improves the quality of  ``Lag. Bound'' as an upper bound in all cases. According to our numerical tests, the optimal value of the inner optimization problem is no greater than zero  for every scenario empirically, and it  becomes farther away from zero with increasing $T$. This happens more often as $\beta$ increases, since larger $\beta$ implies generating a longer horizon $\tau$ with higher probability. Therefore, the improvement of the ``Info. Relaxation''  over ``Lag. Bound'' in absolute value  is more obvious with larger $\beta$.  All the relative gaps are within $2\%$, implying little space of further improvement for both policy and upper bound; comparatively, we may invest more efforts in those problems with smaller $N$ if necessary. In terms of the reduced duality gaps,  the information relaxation bounds improve over the Lagrangian bounds for around  $50\%$ in most cases. This significant improvement shows that the information relaxation approach strengthens the upper bound performance even though the quality of the ``Lagrangian Bound'' has already  been  good.

\subsection{Linear Quadratic Control with Nonconvex linking constraint}
We next consider a finite horizon linear quadratic   control (LQC) problem  with a non-convex linking constraint. We refer the readers to \cite{Haugh:2012} on the information relaxation approach in (unconstrained) finite horizon LQC. Let  $\textbf{x}_t\in \mathcal{X}_t=\mathbb{R}^N$ and $\textbf{a}_t\in\mathcal{A}_t=\mathbb{R}^N$ denote the state and the action at time $t$, respectively. The state equation is described by
\begin{align}
\textbf{x}_{t+1}=A_t\textbf{x}_{t}+B_t \textbf{a}{_t}+\textbf{w}_{t+1}, ~~t=0,\cdots,T-1, \label{LQG_state_equation}
\end{align}
where $A_t,B_t$ are diagonal matrices  for $t=0,\cdots, T-1$, and  $\textbf{w}_{t}'$s are N-dimensional zero-mean  random vectors  with finite second moments. In particular, $cov(\textbf{w}_{t})=\Sigma_t$ is a diagonal matrix for  $t=1,\cdots, T$.  We denote by $\mathbb{F}$ the  natural filtration generated  by $\{\textbf{w}_0,\cdots,\textbf{w}_{T-1}\}$.

The objective is to minimize the expected cost
\begin{align}
U_0(\textbf{x}_0)=\min_{\alpha\in \bar{\mathbb{A}}_{\mathbb{F}}(T)}\mathbb{E}\left[\sum_{t=0}^{T-1}\textbf{a}_t^{\top}\tilde{R}_t \textbf{a}_t+\textbf{x}_T^{\top}Q_T \textbf{x}_T \bigg|\textbf{x}_0\right], \label{LQG_reward}
\end{align}
where each $\tilde{R}_t$ and $Q_T$ are diagonal positive definite matrices, and $\bar{\mathbb{A}}_{\mathbb{F}}(T)$ is the set of non-anticipative policies $\alpha$, where $\alpha$ selects $\textbf{a}=(\textbf{a}_0,\textbf{a}_1,\cdots,\textbf{a}_{T-1})$ over time such that $\textbf{a}_t\in \bar{\mathcal{A}}_t=\{\textbf{a}_t\in \mathbb{R}^N| \tilde{B}(\textbf{a}_t)\triangleq\sum_{n=1}^N (a^n_t)^2 \geq b\}$ with $b\in \mathbb{R}_{+}$ for each $t=0,1,\cdots,T-1$. The system (\ref{LQG_state_equation})-(\ref{LQG_reward}) is weakly-coupled, since $A_t,~B_t,~\Sigma_t,~\tilde{R}_t$, and $Q_T$ are all diagonal matrices and the linking constraint at time $t$ is  $\tilde{B}(\textbf{a}_t)\geq b$. It is simple to verify that the value function $U_0$ is well defined for all $b\geq0$.

Note that the control set $\bar{\mathcal{A}}_t$ is nonconvex, so the optimal policy for (\ref{LQG_reward}) cannot be solved to optimality. Instead we consider a simple heuristic.  At each period $t$ we  compute the one-step greedy policy induced by the value function to the unconstrained problem: we apply such an action if it is already feasible subject to the linking constraint; otherwise, we project it onto the sphere $\partial\bar{\mathcal{A}}_t\triangleq\{\textbf{a}\in \mathbb{R}^N| \tilde{B}(\textbf{a})=b\}$, and use the projection as the action at time $t$. We call this heuristic  ``projection policy''. The performance of this policy provides an upper bound on (\ref{LQG_reward}) (since it  is a minimization problem), which will be referred to as ``Projection Policy'' in Table  \ref{Table:Result2}.

 To derive a lower bound on $U_0$ we first consider the Lagrangian relaxation of  (\ref{LQG_reward}), which  turns out to be an unconstrained LQC problem:
\begin{align*}
J^{\lambda}_0(\textbf{x}_0)&\triangleq \min_{\alpha\in \mathbb{A}_{\mathbb{F}}(T)} \mathbb{E}\left[\sum_{t=0}^{T-1}\textbf{a}_t^{\top}\tilde{R}_t \textbf{a}_t+\textbf{x}_T^{\top}Q_T \textbf{x}_T -\sum_{t=0}^{T-1}\lambda_t\cdot\left[\tilde{B}(\textbf{a}_t)-b\right] \bigg|\textbf{x}_0\right]\\
&= \min_{\alpha\in \mathbb{A}_{\mathbb{F}}(T)} \mathbb{E}\left[\sum_{t=0}^{T-1}\textbf{a}_t^{\top}\left(\tilde{R}_t-\lambda_t\cdot \textbf{I}_N\right) \textbf{a}_t+\textbf{x}_T^{\top}Q_T \textbf{x}_T  \bigg|\textbf{x}_0\right]+\sum_{t=0}^{T-1}\lambda_t^{\top}b,
\end{align*}
where each $\lambda_t$ is a scalar and $\lambda=(\lambda_0,\cdots,\lambda_{T-1})\geq 0$,  and $\textbf{I}_N$ is the N-dimensional identity matrix. Noting that $J^{\lambda}_0(\textbf{x}_0)$  admits a closed form solution  that is  quadratic in $\textbf{x}_0$,  provided that every $\tilde{R}_t-\lambda_t\cdot \textbf{I}_N$ is positive definite:
$$ J^{\lambda}_t(\textbf{x}_0)= \textbf{x}^{\top}_t K_t \textbf{x}_t+ \sum_{s=t}^{T-1} \text{trace}(K_{s+1}\Sigma_{s+1})+\sum_{s=t}^{T-1}\lambda_s\cdot b,~~t=0,\cdots,T.$$
where $K_0$ is obtained by the Riccati equation $K_T=Q_T$, and
$$K_t=A_t'\left(K_{t+1}-K_{t+1}B_t\left(B_t'K_{t+1}B_t+(\tilde{R}_t-\lambda_t\cdot \textbf{I}_N)\right)^{-1}B_t'K_{t+1}\right)A_t,~~t=T-1,\cdots,0.$$

We can  use stochastic subgradient method to derive a tightest Lagrangian bound  on the domain $\mathcal{S}\triangleq\{\lambda\geq0|\tilde{R}_t-\lambda_t\cdot \textbf{I}_N \succ \textbf{0}, t=0,\cdots,T-1\}$. Due to the restricted range, the Lagrangian multiplier $\lambda$ may not be optimal, but $J^{\lambda}_0$ is still  a valid lower bound on $U_0$.

Based on the Lagrangian bounds $\{J^{\lambda}_t\}_{t=1}^T$ we can derive the information relaxation bound through (10) in Appendix \ref{Appendix:FinteHorizon} by choosing $H_t=J^{\lambda}_t(\textbf{x}_t)$, that is,
{\small
\begin{align}\mathbb{E}_0\left[\max_{ \boldsymbol\mu\geq 0}\min_{\textbf{a}\in\mathcal{A}(T)}  \left\{\textbf{x}_T^{\top}Q_T \textbf{x}_T+\sum_{t=0}^{T-1}\textbf{a}_t^{\top}\tilde{R}_t \textbf{a}_t+ \mu_t\cdot(b_t-B_t(\textbf{x}_t,\textbf{a}_t))+ \mathbb{E}[ J^{\lambda}_{t+1}(\textbf{x}_{t+1})|\textbf{x}_t,\textbf{a}_t]- J^{\lambda}_{t+1}(\textbf{x}_{t+1})\right\}\right], \label{LQG_LB}
\end{align}
}where $\mu=(\mu_0,\cdots, \mu_{T-1})$, and
{\small\begin{align*}
\mathbb{E}[ J^{\lambda}_{t+1}(\textbf{x}_{t+1})|\textbf{x}_t,\textbf{a}_t]- J^{\lambda}_{t+1}(\textbf{x}_{t+1}) = -2(A_t\textbf{x}_t+B_t\textbf{a}_t)^{\top}K_{t+1}\textbf{w}_{t+1}-\textbf{w}'_{t+1}K_{t+1}\textbf{w}_{t+1}+  \text{trace}(K_{t+1}\Sigma_{t+1}).
\end{align*}}
Restricting $\mu$ in $\mathcal{S}$, the optimization problem inside the conditional expectation in (\ref{LQG_LB}) is
\begin{align}
\max_{\mu\in\mathcal{S}}\min_{\textbf{a}\in \mathcal{A}(T)} & \bigg\{\textbf{x}_T^{\top}Q_T \textbf{x}_T +\sum_{t=0}^{T-1} \textbf{a}_t^{\top}(\tilde{R}_t-\mu_t\cdot \textbf{I}_N) \textbf{a}_t-2(A_t\textbf{x}_t+B_t\textbf{a}_t)^{\top}K_{t+1}\textbf{w}_{t+1}\notag\\
&-\textbf{w}'_{t+1}K_{t+1}\textbf{w}_{t+1} +  \text{trace}(K_{t+1}\Sigma_{t+1})\bigg\} \label{LQG_inner_opt}
\end{align}
subject to the state dynamics  (\ref{LQG_state_equation}). Then the  minimization problem in (\ref{LQG_inner_opt})  remains a standard deterministic LQ problem, and can be solved efficiently.


{\footnotesize
\begin{table}[htpb]
  \centering \caption{LQ problem with Nonconvex linking constraint } \label{Table:Result2}
  \begin{tabular}{|c|c|c|cc|c|c|cc|cc|} \hline
  \multicolumn{1}{|c} {} &\multicolumn{1}{|c} {} &\multicolumn{1}{|c} {} &\multicolumn{2}{|c|} {Proj. Policy} &\multicolumn{1}{c|} {Unconstrained  } &\multicolumn{1}{c|} {Lag. Bound } &\multicolumn{2}{c|} {Info. Relaxation} &\multicolumn{2}{c|} { Duality Gap }\\ \hline
    { N } &$b$  &$T$ &Value &S.E. &Value   &Value   &Value  &S.E. &1  &2

    \rule[1ex]{0pt}{1em}\\\hline\hline
    $10$   &$5$    &$10$    &$61.4693$     &$0.211$    &$34.7883$     &$59.8636$     &$60.0606$     &$0.0028$   &2.29\% &12.3\% \\
    $20$   &$5$    &$10$    &$ 88.0457$    &$0.242$    &$71.2836$     &$87.2886$     &$87.6371$     &$0.0085$   &0.46\% &46.0\% \\
    $50$   &$5$    &$10$    &$189.1857$    &$0.099$    &$182.7984$    &$188.7715$    &$189.0481$    &$0.0039$   &0.07\% &66.8\%  \\
    $100$  &$5$    &$10$    &$364.2132$    &$0.023$    &$361.7224$    &$364.1160$    &$364.1729$    &$0.0004$   &0.01\% &58.5\%  \\
    $10$   &$10$   &$10$    &$104.6974$    &$0.306$    &$34.7883$     &$103.6067$    &$103.7460$    &$0.0026$   &0.91\% &12.8\%  \\
    $20$   &$10$   &$10$    &$123.4789$    &$0.444$    &$71.2836$     &$ 120.7735$   &$121.5403$    &$0.0090$   &1.57\% &28.3\%  \\
    $50$   &$10$   &$10$    &$209.3848$    &$0.099$    &$182.7984$    &$ 208.8579$   &$209.1757$    &$0.0046$   &0.10\% &60.3\%  \\
    $100$  &$10$   &$10$    &$374.4066$    &$0.193$    &$361.7224$    &$373.7121$    &$374.2227$    &$0.0118$   &0.05\% &73.5\%  \\
    \hline
\end{tabular}
  \end{table}}In our numerical experiments we set $A_t=B_t=\tilde{R}_t=\textbf{I}_N$ for $t=0,\cdots,T-1$, and each diagonal entry of  $Q_T$ is sampled from the uniform distribution on $[1,2]$. We set the initial point $\textbf{x}_0=(1,1,\cdots,1)^{\top}$.  Here is the procedure to get the bounds in Table \ref{Table:Result2}:
\begin{enumerate}
\item[-] ``Proj. Policy'': We generate $10000$ sample paths $\textbf{w}\triangleq(\textbf{w}_0,\cdots,\textbf{w}_{T-1})$ and apply the projection policy to compute the sample cost. To reduce the variance, we use the unconstrained problem as a control variate. The average of the adjusted sample costs provides an upper bound on $U_0$.

\item[-] ``Unconstrained'': The value function to the  problem (\ref{LQG_reward}) without the linking constraint, i.e., $\bar{\mathcal{A}}_t=\mathbb{R}^N$. It can be seen that the ``Unconstrained'' is equal to $J^0$, which is a lower bound on $U_0$.

\item[-] ``Lag. Bound'': we use (stochastic) subgradient method and run $500$ iterations to compute the tightest Lagrangian bound $J^{\lambda^*}_0(\textbf{x}_0)$.  We restrict $\lambda$ in the range $\mathcal{S}'=\{\lambda\geq0|\tilde{R}_t-\lambda_t\cdot \textbf{I}_N\succeq 0.001\cdot\textbf{I}_N,t=0,\cdots,T-1\}\subseteq \mathcal{S}$ (therefore, $\tilde{R}_t-\lambda_t\cdot \textbf{I}_N$ is positive definite) to ease the optimization.  In our numerical experiments the stochastic gradient with respect to  $\lambda$ is very close to zero, which implies that our Lagrangian bound is already near optimal.

\item[-] ``Info. Relaxation'':  We generate another $100$ sample paths of $\textbf{w}$. Based on these sample paths and the Lagrangian bound $J^{\lambda^*}_0$, we compute the relaxed inner optimization problem (\ref{LQG_inner_opt}) (also replace $\mathcal{S}$ by $\mathcal{S}'$) using subgradient method that runs at most $80$ iterations or until the norm of the subgradient is under the tolerance level (we set it to be 0.001).  For most scenarios, this relaxed inner optimization problem can be solved to optimality after around 40 iterations.

\item[-] We also report the duality gaps in two categories:
 \begin{align*}
&\text{Duality Gap 1}=\frac{\text{``Proj.~Policy''~--~``Info.~Relaxation''}}{\text{``Proj.~Policy''}},\\
&\text{Duality Gap 2}=\frac{\text{``Info.~Relaxation''~--~``Lag.~Bound''}}{\text{``Proj.~Policy''~--~``Lag.~Bound''}}.
\end{align*}
\end{enumerate}

Observing  the small gaps between ``Proj. Policy'' and ``Lag. Bound'', it is a little surprising to see the excellent performance of the simple projection policy.  We also note that this simple policy is not trivial by comparing ``Proj. Policy''  to ``Unconstrained'': the weak lower bound of  ``Unconstrained'' indicates that the ``projection'' should occur in some scenarios if not many. The ``Info. Relaxation'' improves the quality of the ``Lag. Bound'', where the duality gaps also behave quite consistently as those in the restless bandits.  The ``Info. Relaxation'' bound shows that the projection policy becomes closer to optimal as $N$ increases. In this example, the linking constraint is a non-decreasing function in the number of subproblems. Therefore, the  linking constraint becomes weaker  as $N$ increases, i.e., the action derived from the unconstrained problem becomes unlikely  to violate the linking constraint.  As we observe, the optimal value to the constrained problem gets closer to the unconstrained one with increasing $N$.

\section{Conclusion}\label{Section:Conclusion}
Lagrangian relaxation and information relaxation are developed to tackle the budget and non-anticipativity constraints that exist universally in general stochastic dynamic programs.  The attraction of studying the interaction of these relaxations particularly in the setting of weakly coupled dynamic programs is  due to the decomposed structure of the Lagrangian bound, as well as the theoretical strong duality guaranteed by the information relaxation. We show that a tighter dual bound, compared with the Lagrangian bound, can be derived  by incorporating it into the information relaxation approach. For large-scale problem, we further develop a computational method to obtain the practical information relaxation bound, which  implies an intermediate relaxation between the Lagrangian and exact information relaxations. The computation of the practical information relaxation bound is easy to implement, and requires little structure of the linking constraints, compared with the approximate linear programming approach that requires designing problem-specific constraint sampling or column generation method.
We may apply this computational method to the case in which both ``easy'' and ``complicated'' linking constraints exist: to balance the complexity and quality of the dual bound, we may choose to only dualize the ``complicated'' constraints in the inner optimization problem.

\bibliographystyle{IEEEtran}
\bibliography{FYe-Bibtex}


%
%
%

\appendices
\section{Approximate Linear Programming Approach}\label{Appendix:ALP}
The approximate linear programming (ALP) method aims to find a good approximation of $V$ within  a parameterized class of functions with a lower-dimensional representation \cite{de2003linear}. In the setting of weakly coupled stochastic dynamic program,  we can set $H(\textbf{x})=\theta+\sum_{n=1}^N H^n(x^n)$, where $\theta$ is a constant and $H^n(\cdot)$ only depends on  $x^n$ for $n=1,\cdots,N$. This approximation scheme is  motivated by the additive form of the Lagrangian bound $J^{\boldsymbol\lambda}$.

Recall that  the set of supersolutions
$$\mathcal{D}^*=\{H\in \mathcal{D}: R(\textbf{x}_0,\textbf{a}_0)+\beta \mathbb{E}[ H(\textbf{x}_{1})|\textbf{x}_0,\textbf{a}_0]\leq H(\textbf{x}_{0}) ~~\text{for all}~ \textbf{x}_0\in  \mathcal{X}  ~\text{and}~ \textbf{a}_0\in \bar{\mathcal{A}}(\textbf{x}_0) \}.$$
Note that each $H^n(\cdot)$ is a mapping from $\mathcal{X}^n$ to $\mathbb{R}$ , which implies that $H(\textbf{x})$ can be parameterized by $1+\sum_{i=1}^n|\mathcal{X}^n|$ variables.
To determine the appropriate parameters, we are seeking a best feasible and additively separable solution  from $\mathcal{D}^*$ via the following linear program \big(with variables $\theta$ and $\{H^{n}(\cdot)\}_{n=1}^N$\big):
\begin{align}
H^{LP}(\upsilon)\triangleq&\min_{\{\theta, H^n(\cdot)\}} \theta+\sum_{n=1}^N\sum_{x^n\in  \mathcal{X}_n} \upsilon_n(x^n) H^n(x^n)   \label{ALP}\\
\text{s.t.} ~~ \theta(1-\beta)+ \sum_{n=1}^N  &H^n(x^n_0)\geq \sum_{n=1}^N R^n(x^n_0,a^n_0)+\beta \sum_{n=1}^N\sum_{x^n_1\in  \mathcal{X}^n}  P_n(x_1^n|x_0^n,a^n_0) H^n(x^n_0), \notag\\
  &~~\text{for all}~~\textbf{x}_0\in\mathcal{X}  ~\text{and}~ \textbf{a}_0\in\bar{\mathcal{A}}(\textbf{x}_0),  \notag
\end{align}
where $\upsilon_n(x^n)$ is the marginal distribution  of $x^n$ from a probability distribution $\upsilon(\cdot)$ on $\mathcal{X}$, and the constraints  are derived from substituting $H(\cdot)$ by $\theta+\sum_{n=1}^N H^n(\cdot)$  in $\mathcal{D}^*$.

We  denote   by $\{\theta^*, H^{LP, n}(\cdot),~n=1,\cdots,N\}$ the optimal solution to (\ref{ALP}), and define
$$H^{LP}(\textbf{x})\triangleq \theta^*+\sum_{n=1}^N H^{LP,n}(x^n).$$

The following lemma shows that the bound derived by the ALP method is tighter than the Lagrangian bound, the proof of which can be found in \cite{adelman2008relaxations}.
\begin{lemma} \label{Proposition:ALP}
\begin{enumerate}
\item[(a)] $H^{LP}\in \mathcal{D}^*$, and $V(\textbf{x})\leq H^{LP}(\textbf{x})$ for all $\textbf{x}\in \mathcal{X}$.
\item[(b)]  $J^{\boldsymbol\lambda}(\textbf{x})=\frac{\boldsymbol\lambda^{\top}\textbf{b}}{1-\beta} \boldsymbol+ \sum_{n=1}^{N}H^{\boldsymbol\lambda,n}(x^n) \in \mathcal{D}^*$, i.e., $\{\frac{\boldsymbol\lambda^{\top} \textbf{b}}{1-\beta}, H^{\boldsymbol\lambda,n}(\cdot),~n=1,\cdots,N\}$  is in the feasible region of the linear program (\ref{ALP}).
\item[(c)] $H^{LP}(\upsilon)\leq J^{\boldsymbol\lambda}(\upsilon)$ for any $\boldsymbol\lambda\in\mathbb{R}^L_{+}$ and probability distribution $\upsilon$.
\end{enumerate}
\end{lemma}

\section{Complements to Section 2} 

\subsection{A formal definition of $\tau$} \label{Appendix:DefinitionTAU}
In this subsection we discuss the augmentation of the probability space $(\Omega, \mathcal{F}, P)$  due to  the introduction of the random time $\tau$.   We can assume that the random variable $\tau$ is associated with another  probability space   $(\hat{\Omega}, \hat{\mathcal{G}},\hat{P})$, where $\tau: \hat{\Omega}\rightarrow \mathbb{N}$, $\hat{\mathcal{G}}$ is the $\sigma$-algebra generated by $\tau$ (i.e., $\sigma(\tau)$), and $\hat{P}(\tau=t)=(1-\beta)\beta^t$ for $t=0,1,2,\cdots$.

The probability space   $(\Omega, \mathcal{F}, P)$ is then augmented  to $(\Omega\times \hat{\Omega}, \mathcal{F}\otimes\sigma(\tau), \mathbb{P})$, where $\mathcal{F}\otimes\sigma(\tau)$ is the product  $\sigma$-algebra of  $\mathcal{F}$ and $\sigma(\tau)$, and $\mathbb{P}$ is the product measure of $P$ and $\hat{P}$, i.e., $\mathbb{P}(A\times [t,\infty))=P(A)\times \hat{P}(\tau\geq t)=P(A)\times \beta^{t}$ with $A\in \mathcal{F}$. We clarify this (straightforward) augmentation  is because we can use the pair $(\omega,\tau)$ to denote the uncertainty in the conditional expectation in $\mathcal{L}H$  without confusion,  though to save notations we use $P$ to denote  $\mathbb{P}$.




\subsection{$\max_{\textbf{a}\in \mathcal{\bar{A}}(\tau)}\{I_H(\textbf{a},\omega,\tau)\}$ has finite mean and variance}\label{Appendix:welldefine}
Let $\mathcal{I}(\omega,\tau)=\max_{\textbf{a}\in \mathcal{\bar{A}}(\tau)}\{I_H(\textbf{a},\omega,\tau)\}$. Then $\mathcal{L}H(\textbf{x}_0)\triangleq H(\textbf{x}_{0})+\mathbb{E}_0[\mathcal{I}(\omega,\tau)]$. Since $R$ and $H$ are both bounded, we can assume for all $(\textbf{x}_t,\textbf{a}_t)\in \mathcal{X}\times\mathcal{ A}$, $t=0,1,2,\cdots,$
$$|R(\textbf{x}_t,\textbf{a}_t)+\beta \mathbb{E}[ H(\textbf{x}_{t+1})|\textbf{x}_t,\textbf{a}_t]- H(\textbf{x}_{t})|\leq C$$
for some $C>0$. Therefore,  $|\mathcal{I}(\omega,\tau)|\leq(\tau+1)C$ for any $\omega\in \Omega$, which implies
\begin{equation}
|\mathbb{E}_0\left[\mathcal{I}(\omega,\tau)\right]|\leq\mathbb{E}_0\left[\mathbb{E}\big[|\mathcal{I}(\omega,\tau)|\big|\tau\big]\right]\leq \sum_{\tau=0}^{\infty}(1-\beta)\beta^{\tau}(\tau+1)C=\frac{C}{1-\beta}<\infty, \label{finite_mean}
\end{equation}
and $\text{Var}[\mathcal{I}(\omega,\tau)|\tau]\leq \mathbb{E}\left[\mathcal{I}^2(\omega,\tau)|\tau\right]\leq (\tau+1)^2C^2. $ The inequality (\ref{finite_mean}) indicates that $\mathcal{I}(\omega,\tau)$ has finite mean.

We  note that $\text{Var}[\mathcal{I}(\omega,\tau)]=\mathbb{E}\left[\text{Var}\big[\mathcal{I}(\omega,\tau)\big|\tau\big]\right]+\text{Var}\left[\mathbb{E}\big[\mathcal{I}(\omega,\tau)\big|\tau\big]\right].$

It can be seen that $\mathbb{E}[\text{Var}[\mathcal{I}(\omega,\tau)|\tau]]\leq \sum_{\tau=0}^{\infty}(1-\beta)\beta^{\tau}(\tau+1)^2C^2=\frac{1+\beta}{(1-\beta)^2}C^2<\infty,$
and
$$\text{Var}[\mathbb{E}[\mathcal{I}(\omega,\tau)|\tau]]\leq \mathbb{E}[(\mathbb{E}[\mathcal{I}(\omega,\tau)|\tau])^2]\leq \mathbb{E}[(\tau+1)^2C^2]=\sum_{\tau=0}^{\infty}(1-\beta)\beta^{\tau}(\tau+1)^2C^2=\frac{1+\beta}{(1-\beta)^2}C^2< \infty.$$
Hence, we conclude that $\mathcal{I}(\omega,\tau)$ has finite variance.

\section{Information Relaxation Improves the Lagrangian Bound: An Example} \label{Section:Information_Bound:Example}
We consider the restless bandit-like problem with $N=1$ as proposed in Section 3.3 of \cite{adelman2008relaxations}: the state space contains three states, i.e., $\mathcal{X}=\{0,1,2\}$, and for each state $\textbf{x}\in\mathcal{X}$ the control space is $\mathcal{A}(\textbf{x})=\{0,1\}$. The corresponding reward $R(\textbf{x},\textbf{a})$, weight $\textbf{B}(\textbf{x},\textbf{a})$, and transition probability $P(\textbf{x}_{t+1}|\textbf{x}_t,\textbf{a}_t)$ are listed  in Table  \ref{Table:output1}, in which $l>0$ and $c>1$ are positive constants.  Note that states ``1'' and ``2'' are absorbing states regardless of the control applied; however, the state ``0'' may transit to either ``1'' or ``2'' depending on the control chosen. The linking constraint is $\textbf{B}(\textbf{x},\textbf{a})\leq 1$. Therefore, $\mathcal{\bar{A}}(0)=\mathcal{\bar{A}}(2)=\{0,1\}$ and $\mathcal{\bar{A}}(1)=\{0\}$.

\begin{remark}
In Table 1 of \cite{adelman2008relaxations}, $\textbf{B}(2,0)=\epsilon>0$. All results therein  are also true for $\epsilon=0$.
\end{remark}
\begin{table}[htpb]
  \centering \caption{One-subproblem with $b=1$ and $\beta\in (\frac{1}{2},1)$ } \label{Table:output1}
  \begin{tabular}{|c|c|c|c|c|} \hline
    { State } &Control  &Reward &Weight &Transition

    \rule[1ex]{0pt}{1em}\\\hline\hline

    $0$      &$0$     &$R(0,0)=0$       &$\textbf{B}(0,0)=0$     &$P(2|0,0)=1$    \\
    $0$      &$1$     &$R(0,1)=0$       &$\textbf{B}(0,1)=0$     &$P(1|0,1)=1$    \\
    \hline
    $1$      &$0$     &$R(1,0)=0$       &$\textbf{B}(1,0)=0$     &$P(1|1,0)=1$     \\
    $1$      &$1$     &$R(1,0)=c(2+l)$  &$\textbf{B}(1,1)=2$     &$P(1|1,1)=1$     \\
    \hline
    $2$      &$0$     &$R(2,0)=0$       &$\textbf{B}(2,0)=0$     &$P(2|2,0)=1$       \\
    $2$      &$1$     &$R(2,1)=c$       &$\textbf{B}(2,1)=0$     &$P(2|2,1)=1$       \\
\hline
\end{tabular}
\end{table}
The exact value function  is $ V(0)=\frac{c\beta}{1-\beta},~ V(1)=0, ~\text{and} ~V(2)=\frac{c}{1-\beta}.$
The optimal stationary policy is $\alpha=(\alpha^*_{\delta},\alpha^*_{\delta},\cdots)$, where  $\alpha^*_{\delta}(0)=\alpha^*_{\delta}(1)=0$ and $\alpha^*_{\delta}(2)=1$.  

The Lagrangian relaxation yields $J^{\lambda}(\textbf{x})= \frac{\lambda}{1-\beta}+H^{\lambda}(\textbf{x})$ for $\textbf{x}=0,1,2$. According to \cite{adelman2008relaxations}, the optimal Lagrangian multiplier  is
$$\lambda^{*}=\arg \min_{\lambda\geq 0} H^{\lambda}(\upsilon) = c+cl/2,$$
which implies  $H^{\lambda^*}(0)=0$, $H^{\lambda^*}(1)=0$, and  $H^{\lambda^*}(2)=0$.
Therefore,
$$J^{\lambda^*}(0)= \frac{\lambda^*}{1-\beta},~ J^{\lambda^*}(1)= \frac{\lambda^*}{1-\beta}, ~\text{and}~ J^{\lambda^*}(2)= \frac{\lambda^*}{1-\beta}.$$
Note that $J^{\lambda^*}(\cdot)$ is unbounded on $\mathcal{X}$ as $l\rightarrow \infty,$ though the exact values $V(\cdot)$ is constant with respect to $l$.

By applying the information relaxation approach with $H=J^{\lambda^{*}}$,
\begin{align*}
\mathcal{L}J^{\lambda^*}(\textbf{x}_0)= &H(\textbf{x}_{0})+\mathbb{E}_0\left[\max_{ \textbf{a}\in \mathcal{\bar{A}}(\tau)}\left\{\sum_{t=0}^{\tau}\big(R( \textbf{x}_t,\textbf{a}_t)+\beta \mathbb{E}[ H(\textbf{x}_{t+1})|\textbf{x}_t,\textbf{a}_t]- H(\textbf{x}_{t})\big)\right\}\right]\\
=&J^{\lambda^{*}}(\textbf{x}_{0})+\sum_{T=0}^{\infty}(1-\beta)\beta^{T} \cdot  \mathbb{E}_0\left[\max_{\textbf{a}\in \mathcal{\bar{A}}(\tau)} \left\{\sum_{t=0}^{T}\big(R(\textbf{x}_t,\textbf{a}_t)+\beta \mathbb{E}[ J^{\lambda^*}(\textbf{x}_{t+1})|\textbf{x}_t,\textbf{a}_t]- J^{\lambda^*}(\textbf{x}_{t})\big)\right\}\right].
\end{align*}

Note that  $J^{\lambda^*}(\textbf{x}_0)- (\textbf{R}(\textbf{x}_0,\textbf{a}_0)+\beta \mathbb{E}[ J^{\lambda^*}( \textbf{x}_1)|\textbf{x}_0,\textbf{a}] ) \geq  \lambda^*-c$ for all  $\textbf{x}_0\in \mathcal{X}=\{1,2,3\}$  and $\textbf{a}_0\in \bar{\mathcal{A}}(\textbf{x}_0)$.  According to Theorem 2(c), $J^{\lambda^*}(\textbf{x})- \mathcal{L}J^{\lambda^*}(\textbf{x})\geq \frac{\lambda^*-c}{1-\beta}$, which implies  $\mathcal{L}J^{\lambda^*}(\textbf{x})\leq \frac{c}{1-\beta}$  for $\textbf{x}\in \mathcal{X}$. This bound remains constant with respect to $l$ and it has already been tight as an upper bound on $V(2)$.

We can show that the exact computation of the information relaxation bound also leads to a tight upper bound on $V(0)$. Starting at $\textbf{x}_0=0$ and for each $T\in \mathbb{N}$ and $\omega\in \Omega$,
\begin{align*}
\max_{\textbf{a}\in\mathcal{A}(T)}&\left\{\sum_{t=0}^{T}\big(R(\textbf{x}_t,\textbf{a}_t)+\beta \mathbb{E}[ J^{\lambda^*}(\textbf{x}_{t+1})|\textbf{x}_t,\textbf{a}_t]-
J^{\lambda^*}(\textbf{x}_{t})\big)\right\}\\
=\max &\bigg \{R(0,0)+\beta \mathbb{E}[ J^{\lambda^*}(\textbf{x}_{1})|0,0]- J^{\lambda^*}(0)+\sum_{t=1}^{T} \max_{\textbf{a}_t\in \bar{\mathcal{A}}(\textbf{x}_t)} \left\{ R(
\textbf{x}_t, \textbf{a}_t)+\beta \mathbb{E}[ J^{\lambda^*}(\textbf{x}_{t+1})|\textbf{x}_t,\textbf{a}_t]- J^{\lambda^*}(\textbf{x}_{t})\right\}, \\
&R(0,1)+\beta \mathbb{E}[ J^{\lambda^*}(x_{1})|0,1]- J^{\lambda^*}(0)+\sum_{t=1}^{T}\max_{\textbf{a}_t\in \bar{\mathcal{A}}(\textbf{x}_t)} \big\{ R(\textbf{x}_t,\textbf{a}_t)+\beta
\mathbb{E}[ J^{\lambda^*}(\textbf{x}_{t+1})|\textbf{x}_t,\textbf{a}_t]- J^{\lambda^*}(\textbf{x}_{t})\big\} \bigg\}\\
=\max &\bigg\{R( 0, 0)+\beta J^{\lambda^*}(2)- J^{\lambda^*}(0)+\sum_{t=1}^{T}\max_{\textbf{a}_t\in \bar{\mathcal{A}}(2)} \left\{R(2, \textbf{a}_t)- (1-\beta) J^{\lambda^*}( 2)\right\},
\\
 &R(0,1)+\beta J^{\lambda^*}(1)- J^{\lambda*}(0)+\sum_{t=1}^{T}\max_{\textbf{a}_t\in \bar{\mathcal{A}}(1)} \left\{R( 1, \textbf{a}_t)- (1-\beta) J^{\lambda^*}(1)\right\} \bigg\}\\
=\max &\left\{0+\beta \frac{\lambda^*}{1-\beta}-\frac{\lambda^*}{1-\beta} +(c-\lambda^*)T,~0+\beta \frac{\lambda^*}{1-\beta}-\frac{\lambda^*}{1-\beta}+(0-\lambda^*)T\right\}\\
=-\lambda^*&+(c-\lambda^*)T,
\end{align*}
where the first equality holds since staring at $\textbf{x}_0=0$, the control $\textbf{a}_0=0$  leads to  $\textbf{x}_1=2$ (respectively, $\textbf{a}_0=1$  leads to  $\textbf{x}_1=1$) with probability $1$, and hence determine all the subsequent states $\textbf{x}_2$, $\textbf{x}_3$, $\cdots$, since  $\textbf{x}=1$ and $2$ are absorbing states. Consequently, the deterministic dynamic program with time horizon $T$ can be decomposed as the summation of $T$ sub-problems.  The last equality holds as the first term dominates the second, meaning that $\textbf{a}_0=0$ and $\textbf{a}_1=1$ for $t\geq1$ is the solution to the inner optimization problem for all the scenarios $\omega\in \Omega$. Since  $J^{\lambda^*}(0)=\frac{\lambda^*}{1-\beta}$, then
\begin{align*}
\mathcal{L}J^{\lambda^*}(0)=&\frac{\lambda^*}{1-\beta}+ \mathbb{E}_0[-\lambda^* +(c-\lambda^*)\tau]=\frac{\lambda^*}{1-\beta}+ \sum_{\tau=0}^{\infty}(1-\beta)\beta^{\tau} [-\lambda^* +(c-\lambda^*)\tau]=\frac{ c\beta}{1-\beta}.
\end{align*}
Hence, $\mathcal{L}J^{\lambda^*}(0)=V(0)$.

Note the solution to the inner optimization problem of any horizon $T$  is of the form $\textbf{a}=(0,1,1,1,\cdots)$, and the resulting trajectory of the state $(\textbf{x}_0,\textbf{x}_1,\textbf{x}_2,\textbf{x}_3,\cdots)$ is  $(0,2,2,2,\cdots)$ for all $\omega\in\Omega$ and $T\geq1$. This confirms the conditions in Theorem 3, as the optimal policy to the original problem is  $\alpha^*_{\delta}(0)=0$, $\alpha^*_{\delta}(1)=0$,  and $\alpha^*_{\delta}(2)=1$. In particular, $\alpha^*_{\delta}$ is exactly the greedy policy induced by the Lagrangian bound $J^{\lambda^*}$.

\section{Proof of Theorem 4}\label{Appendix:ProofTheorem4}
\begin{IEEEproof}
Given $\alpha'\in \mathbb{\bar{A}}_{\mathbb{F}}$ and $\textbf{x}_0\in\mathcal{X}$,
\begin{align}
V(\textbf{x}_0;\alpha')=&\mathbb{E}_0\left[\sum_{t=0}^{\infty}\beta^t R(\textbf{x}_t,\alpha'_{\delta}(\textbf{x}_t))\right] \notag\\
=&\mathbb{E}_0\left[\sum_{t=0}^{\infty}\beta^t R(\textbf{x}_t,\alpha'_{\delta}(\textbf{x}_t))-  \beta\cdot \sum_{t=0}^{\infty}\beta^t \Delta H_{t+1}(\alpha',\omega)\right]  \notag \\
=&H(\textbf{x}_{0})+\mathbb{E}_0\left[\sum_{t=0}^{\tau}\big(R(\textbf{x}_t,\alpha'_{\delta}(\textbf{x}_t))+\beta \mathbb{E}[ H(\textbf{x}_{t+1})|\textbf{x}_t,\alpha'_{\delta}(\textbf{x}_t)]- H(\textbf{x}_{t})\big)\right] \notag\\
\leq &H(\textbf{x}_{0})+\mathbb{E}_0\left[\max_{\textbf{a}\in \mathcal{\bar{A}}(T)} \left\{\sum_{t=0}^{\tau}\big(R(\textbf{x}_t,\textbf{a}_t)+\beta \mathbb{E}[
H(\textbf{x}_{t+1})|\textbf{x}_t,\textbf{a}_t]- H(\textbf{x}_{t})\big)\right\}\right],  \label{context_complementary22}\\
=&H(\textbf{x}_{0})+\sum_{T=0}^{\infty}P(\tau=T)\cdot\max_{\textbf{a}\in \mathcal{\bar{A}}(T)} \left\{\sum_{t=0}^{T}\big(R(\textbf{x}_t,\textbf{a}_t)+\beta \mathbb{E}[
H(\textbf{x}_{t+1})|\textbf{x}_t,\textbf{a}_t]- H(\textbf{x}_{t})\big)\right\} \notag\\
=&\mathcal{L}H(\textbf{x}_0).  \notag
\end{align}

To show necessity,  $V(\textbf{x};\alpha')=\mathcal{L}H(\textbf{x})$ means that the inequality (\ref{context_complementary22}) is an equality; by observing that $P(\tau=T)>0$ for every $T\in \mathbb{N}$, the equality(\ref{context_complementary21}) should hold for $\omega\in \Omega$ almost surely, $T=0,1,2,\cdots.$

%
%
The sufficiency is straightforward, since the condition (\ref{context_complementary21}) holds for $\omega\in \Omega$ almost surely and $T\in \mathbb{N}$ implies that (\ref{context_complementary22}) is an equality, and thus $V(\textbf{x}_0;\alpha')=\mathcal{L}H(\textbf{x}_0)$.
\end{IEEEproof}

\section{Proof of Theorem 5}\label{Appendix:Separable_DulaityGap}
To prove Theorem 5, we use the result of Lagrangian duality gap on deterministic separable problem. Consider a separable problem
\begin{align}
\max_{\textbf{a}\in \bar{\mathcal{A}}}\sum_{n=1}^{N}f^n(a^n), \label{det_sep}
\end{align}
where $\bar{\mathcal{A}}=\{\textbf{a}\triangleq (a^1,\cdots,a^N)\in \mathcal{A}^1\times\cdots \times \mathcal{A}^N|~\sum_{n=1}^{N}\textbf{h}^n(a^n)\leq \textbf{q}\}$  with $\textbf{q}\in \mathbb{R}^{\tilde{L}}$.

We then define the Lagrangian dual of (\ref{det_sep}):
$$\min_{\mu\geq 0}~d(\mu)\triangleq\sum_{n=1}^N \max_{a^n\in  \mathcal{A}^n }\{f^n(a^n)-\mu^{\top}\textbf{h}^n(a^n)\}+\mu^{\top}\textbf{q}.$$

\begin{lemma}[Proposition 5.26 in \cite{bertsekas1982constrained}] \label{Lemma:Separable_DulaityGap}
 Suppose the following assumptions hold.
 \begin{enumerate}
 \item[] Assumption 1: $\bar{\mathcal{A}}\neq \emptyset$.
 \item[] Assumption 2: for each $n=1,\cdots,N$,  $\{a^n, \textbf{h}^n(a^n), f^n(a^n)|a^n\in\mathcal{A}^n \}$  is compact.
 \item[] Assumption 3: for each $n=1,\cdots,N$, given any vector $\tilde{a}^n\in conv(\mathcal{A}^n)$, there exists $a^n\in \mathcal{A}^n$ such that
 $$\textbf{h}^n(a^n)\leq (\check{\text{cl}}~\textbf{h}^n)(\tilde{a}^n).$$
 \end{enumerate}
 Then
$$\min_{\mu\geq 0} d(\mu)-\max_{\textbf{a}\in \bar{\mathcal{A}}}\sum_{n=1}^{N}f^n(a^n)\leq (\tilde{L}+1)\max_{n=1,\cdots,N}\rho_n,$$
 where  $\rho_n=\sup_{a^n\in conv(\mathcal{A}^n)}\left\{ \widetilde{f}^n(a^n)-(\check{\text{cl}}~f^n)(a^n) \right\}$.
\end{lemma}

The proof of Theorem 5 uses the  following lemma, which is a corollary of Lemma \ref{Lemma:Separable_DulaityGap}.
\begin{lemma}\label{Lemma:Inner_Separable_DulaityGap}
Suppose that $H$ is  of the additively separable form  $H(\textbf{x})=\theta+\sum_{n=1}^N H^n(x^n)$, and Assumptions 1-3 in Section 3.2 hold for $\omega\in \Omega$ and $T\in \mathbb{N}$. Then
$$\min_{ \boldsymbol\mu\geq 0} \max_{\textbf{a}\in\mathcal{A}(T)}I_{H}(\textbf{a},\omega,T; \boldsymbol\mu)-\max_{\textbf{a}\in \mathcal{\bar{A}}(T)} I_{H}(\textbf{a},\omega,T)\leq (1+ L(T+1) )\max_{n=1,\cdots, N} \gamma^n,$$
where
$$\gamma^n=  \sup_{\tilde{\textbf{a}}^n\in conv(\mathcal{A}^n(T))}\left\{\widetilde{I^n_{H_n}}(\tilde{\textbf{a}}^{n},\omega,T;0)- (\check{\text{cl}}~I^n_{H_n})(\tilde{\textbf{a}}^{n},\omega,T;0)\right\},$$
$\check{\text{cl}}~I^n_{H}$ is the convex closure of $I^n_{H_n}$, and $\widetilde{I^n_{H_n}}$ is defined as
$$\widetilde{I^n_{H^n}}(\tilde{\textbf{a}}^{n},\omega,T;0)=\inf_{\textbf{a}^n\in \mathcal{A}^n(T)}\left\{I^n_{H^n}(\textbf{a}^{n},\omega,T;0)|\textbf{B}^n_t(\textbf{a}^n, \omega)\leq  (\check{\text{cl}}~\textbf{B}^n_t)(\tilde{\textbf{a}}^n, \omega), ~t=0,\cdots,T\right\}.$$
\end{lemma}
\begin{remark}
Note that $\widetilde{I^n_{H^n}}(\tilde{\textbf{a}}^{n},\omega,T;0)$ is well-defined according to Assumption 3 in Section 3.2.
\end{remark}
\begin{IEEEproof}
Lemma \ref{Lemma:Inner_Separable_DulaityGap} directly follows from Lemma \ref{Lemma:Separable_DulaityGap} by setting $f^n=I^n_{H^n}$, $\textbf{h}^n=(\textbf{B}^n_0,\cdots,\textbf{B}^n_T)$, $\textbf{q}=(\textbf{b},\cdots,\textbf{b})\in \mathbb{R}^{\tilde{L}}$ with $\tilde{L}=L\times(T+1)$, and the decision variable $a^n=\textbf{a}^n\in \mathcal{A}^n(T)$.
\end{IEEEproof}

\begin{thm7}\label{Theorem:Bound_Gap}
Suppose that $H$ is  of the additively separable form  $H(\textbf{x})=\theta+\sum_{n=1}^N H^n(x^n)$, and Assumptions (1)-(3) hold for every $\omega\in \Omega$ and $T\in \mathbb{N}$. Then for all $\textbf{x}\in \mathcal{X}$,
  \begin{equation}
  \mathcal{L}^\circ H(\textbf{x})-\mathcal{L}H(\textbf{x})\leq \frac{(L-1)\beta+L+1}{(1-\beta)^2}  \max_{n=1,\cdots, N} \Gamma^n, \label{Gap_Bound}
  \end{equation}
where
\begin{align*}
\Gamma^n=&\sup_{x^n_0\in \mathcal{X}^n, a^n_0 \in \mathcal{A}^n(x^n_0)}\{ R^n(x^n_0,a^n_0)+\beta \mathbb{E}[ H^n(x^n_{1})|x^n_0,a^n_0]- H^n(x^n_{0}) \}\\
&-\inf_{x^n_0\in \mathcal{X}^n, a^n_0 \in \mathcal{A}^n(x^n_0)}\{ R^n(x^n_0,a^n_0)+\beta \mathbb{E}[ H^n(x^n_{1})|x^n_0,a^n_0]- H^n(x^n_{0})\}.
\end{align*}
\end{thm7}
\begin{IEEEproof}
According to Lemma \ref{Lemma:Inner_Separable_DulaityGap}, we have for fixed $\omega\in\Omega$ and $\tau=T$,
$$\min_{\mu\geq 0} \max_{\textbf{a}\in\mathcal{A}(T)}I_{H}(\textbf{a},\omega,T;\mu)-\max_{\textbf{a}\in\mathcal{\bar{A}}(T)} I_{H}(\textbf{a},\omega,T)\leq (1+ L(T+1) )\max_{n=1,\cdots, N} \gamma^n,$$
 where
\begin{align*}
\gamma^n \leq &\sup_{\textbf{a}^{n}\in\mathcal{A}^n(\tau)}\{I^n_{H^n}(\textbf{a}^{n},\omega,T;0)\}-  \inf_{\textbf{a}^{n}\in \mathcal{A}^n(\tau)}\{ I^n_{H^n}(\textbf{a}^{n},\omega,T;0)\}\\
\leq &(T+1) \sup_{x^n_0\in \mathcal{X}^n, a^n_0 \in \mathcal{A}^n(x^n_0)}\{ R^n(x^n_0,a^n_0)+\beta \mathbb{E}[ H^n(x^n_{1})|x^n_0,a^n_0]- H^n(x^n_{0}) \}\\
&- (T+1) \inf_{x^n_0\in \mathcal{X}^n, a^n_0 \in \mathcal{A}^n(x^n_0)}\{ R^n(x^n_0,a^n_0)+\beta \mathbb{E}[ H^n(x^n_{1})|x^n_0,a^n_0]- H^n(x^n_{0})\}\\
=&(T+1) \Gamma^n,
\end{align*}
where the first inequality is due to the definitions of $\widetilde{I^n_{H^n}}$ and $\check{\text{cl}}~I^n_{H^n}$, and the second inequality holds independent of $\omega$.
 It is straightforward to see
\begin{align*}
\mathcal{L}^\circ H(x)-\mathcal{L}H(x) &= \mathbb{E}\left[\min_{ \boldsymbol\mu\geq 0} \max_{\textbf{a}\in \mathcal{A}(\tau)}\{I_{H}(\textbf{a},\omega,\tau; \boldsymbol\mu)\}-\max_{\textbf{a}\in \mathcal{\bar{A}}(\tau)} \{I_{H}(\textbf{a},\omega,\tau)\}\right]\\
&\leq \mathbb{E}\left[\mathbb{E}\bigg[(1+L(\tau+1))(\tau+1) \max_{n=1,\cdots, N} \Gamma^n \bigg|\tau\bigg]\right].
\end{align*}

Then we can obtain  (\ref{Gap_Bound}), since
\begin{align*}
\mathbb{E}\left[\mathbb{E}\left[(1+L(\tau+1))(\tau+1) \max_{n=1,\cdots, N} \Gamma^n \bigg|\tau\right]\right]&=  \max_{n=1,\cdots, N} \Gamma^n \cdot \mathbb{E}\left[(1+L(\tau+1))(\tau+1)\right]\\
&= \frac{(L-1)\beta+L+1}{(1-\beta)^2} \max_{n=1,\cdots, N} \Gamma^n. 
\end{align*}
\end{IEEEproof}

\section{Finite horizon case}\label{Appendix:FinteHorizon}
In this section we  consider the finite-horizon weakly coupled dynamic program, which is the same as infinite-horizon case except that
\begin{enumerate}
 \item The time is indexed by $t=0,\cdots,T.$
 \item The transition probability can be time-varying.
 \item The linking constraint can be time-varying, and the feasible control set at time $t$ is
 \begin{align*}
 \bar{\mathcal{A}}_t(\textbf{x}_t)=\{\textbf{a}=(a^1_t,\cdots,a^N_t)\in \mathcal{A}_t(\textbf{x}_t) :~ \textbf{B}_t(\textbf{x}_t,\textbf{a}_t)\triangleq\sum_{n=1}^N \textbf{B}^{n}_t(x^n_t,a^n_t)\leq \textbf{b}_t \},
 \end{align*}
 where each $\textbf{b}_t\in\mathbb{R}^L$ for $t=0,\cdots,T$.
 \item The intermediate rewards denoted by $R_t(x_t,a_t)=\sum_{n=1}^N R^n_t(x^n_t,a^n_t)$ can also be time-varying.
\end{enumerate}

The objective of the decision maker is to maximize the expected rewards given $\textbf{x}_0\in\mathcal{X}$,
\begin{align}
U_0(\textbf{x}_0)=\max_{\alpha\in \bar{\mathbb{A}}_{\mathbb{F}}(T)}U_0(\textbf{x}_0; \alpha), \label{context_finite_value_function}
\end{align}
where \begin{align*}
U_0(\textbf{x}_0; \alpha)=\mathbb{E}\left[\sum_{t=0}^{T}R_t(\textbf{x}_t,\textbf{a}_t)\bigg|\textbf{x}_0\right], 
\end{align*}
and $\bar{\mathbb{A}}_{\mathbb{F}}(T)$ is the set of non-anticipative policies $\alpha$  that selects $\textbf{a}_t\in \bar{\mathcal{A}}_t(\textbf{x}_t)$ for each $t=0,1,\cdots,T$. Then $U_0$ can be solved via the dynamic programming:
\begin{align*}
U_{T+1}(\textbf{x}_{T+1})&=0;\\
U_t(\textbf{x}_t)&=\max_{\textbf{a}_{t}\in \bar{\mathcal{A}}_t(\textbf{x}_t)}\left\{R_{t}(\textbf{x}_{t},\textbf{a}_{t})+\mathbb{E}[U_{t+1}(\textbf{x}_{t+1})|\textbf{x}_{t},\textbf{a}_{t}]\right\}.
\end{align*}

\subsubsection{Lagrangian Relaxation}
Let $\mathbb{A}_{\mathbb{F}}(T)=\{\alpha\in \mathbb{A}(T)|~\alpha~\text{is non-anticipative}\}.$
By  dualizing the linking constraint with Lagrangian multipliers $\boldsymbol\lambda=(\boldsymbol\lambda_0,\cdots,\boldsymbol\lambda_T)\geq0$ with each $\boldsymbol\lambda_t\in\mathbb{R}^L_{+}$,  we define for $\textbf{x}_0\in \mathcal{X}$,
\begin{align}
J^{\boldsymbol\lambda}_0(\textbf{x}_0)\triangleq \max_{\alpha\in \mathbb{A}_{\mathbb{F}}(T)}J^{\boldsymbol\lambda}_0(\textbf{x}_0;\alpha) \label{context_Lagrangian1},
\end{align}
where
$$J^{\boldsymbol\lambda}_0(\textbf{x}_0;\alpha)\triangleq\mathbb{E}\left[\sum_{t=0}^{T} R_t( \textbf{x}_t, \textbf{a}_t)+\boldsymbol\lambda_t^{\top}\left[\textbf{b}_t-\textbf{B}_t(\textbf{x}_t,\textbf{a}_t)\right] \bigg|\textbf{x}_0\right],$$
and $\mathbb{A}_{\mathbb{F}}(T)$ is the set of non-anticipative policies $\alpha$  that selects $\textbf{a}_t\in \mathcal{A}_t(\textbf{x}_t)$ for each $t=0,\cdots,T$. Then $J_0^{\boldsymbol\lambda}$ can be solved via the dynamic programming equations:
\begin{align}
J_{T+1}^{\boldsymbol\lambda}(\textbf{x}_{T+1})&=0; \notag\\
J_t^{\boldsymbol\lambda}(\textbf{x}_t)&=\max_{\textbf{a}_{t}\in \mathcal{A}_t(\textbf{x}_t)}\left\{R_t(\textbf{x}_{t},\textbf{a}_{t})+\boldsymbol\lambda_t^{\top}[\textbf{b}_t-\textbf{B}_t(\textbf{x}_t,\textbf{a}_t)]+\mathbb{E}[J_t^{\boldsymbol\lambda}(\textbf{x}_{t+1})|\textbf{x}_{t},\textbf{a}_{t}]\right\}. \label{finitehorizon_Lag_DP}
\end{align}

Similar to the infinite-horizon case, the solution to (\ref{context_Lagrangian1}) can be solved by decomposing (\ref{finitehorizon_Lag_DP}) into $N$ dynamic programs of lower dimensions:
$$J^{\boldsymbol\lambda}_{0}(\textbf{x}_{0};\alpha)=\sum_{t=0}^{T} \boldsymbol\lambda_t^{\top}\textbf{b}_t+ \mathbb{E}\left[\sum_{t=0}^{T} R_t( \textbf{x}_t, \textbf{a}_t)-\boldsymbol\lambda_t^{\top}\textbf{B}_t(\textbf{x}_t,\textbf{a}_t) \bigg|\textbf{x}_0\right]
=\sum_{t=0}^{T} \boldsymbol\lambda_t^{\top}\textbf{b}_t+\sum_{n=1}^N H^{\boldsymbol\lambda,n}_{0}(x^n_{0}),$$
where
\begin{align*}
H^{\boldsymbol\lambda,n}_{T+1}(x^n_{T+1})&=0,\\
H^{ \boldsymbol\lambda,n}_t(x^n_t)&=\max_{a^n_{t}\in \mathcal{A}^n_t(x^n_t)}\left\{R_t^n(x^n_{t},a^n_{t})-\boldsymbol\lambda_t^{\top}\boldsymbol B^n_t(x^n_t,a^n_t)+\mathbb{E}[H^{ \boldsymbol\lambda,n}_{t+1}(x^n_{t+1})|x^n_{t},a^n_{t}]\right\}.
\end{align*}

\subsubsection{Information Relaxation}
We define the space of a sequence of  functions $H=(H_0,\cdots,H_{T+1})$:
$$\mathcal{D}_T\triangleq\{H=(H_0,\cdots,H_{T+1})| H_t:\mathcal{X}\rightarrow \mathbb{R} ~\text{for}~ t=0,\cdots,T+1, and~H_{T+1}(\cdot)\equiv 0\}.$$
Given $H\in\mathcal{D}_T$, we define
\begin{align*}
\mathcal{L}_{T}H(\textbf{x}_0) \triangleq & \mathbb{E}_0\left[\max_{\textbf{a}\in \mathcal{A}(T)} \left\{\sum_{t=0}^{T}\big(R_t(\textbf{x}_t,\textbf{a}_t)+ \mathbb{E}[ H_{t+1}(\textbf{x}_{t+1})|\textbf{x}_t,\textbf{a}_t]- H_{t+1}(\textbf{x}_{t+1})\big)\right\}\right]\\
=&H_0(\textbf{x}_{0})+\mathbb{E}_0\left[\max_{\textbf{a}\in \mathbb{A}(T)} \left\{I_{H}(\textbf{a},\omega,T)\right\}\right],
\end{align*}
where we redefine $\textbf{a}\triangleq(\textbf{a}_0,\cdots,\textbf{a}_T)$, and
$$I_{H}(\textbf{a},\omega,T)\triangleq\sum_{t=0}^{T}\big(R_t(\textbf{x}_t,\textbf{a}_t)+ \mathbb{E}[ H_{t+1}(\textbf{x}_{t+1})|\textbf{x}_t,\textbf{a}_t]- H_t(\textbf{x}_{t})\big).$$

\textbf{Practical Information Relaxation Bound}  We further assume for  each $t=0,\cdots, T$, the function $H_{t}$ is of the additively separable form
$$H_{t}(\textbf{x}_t)=\theta_t+\sum_{n=1}^N H^n_{t}(x_t^n),$$
where $\theta_t\in \mathbb{R}$ and $H^n_{t}:\mathcal{X}^n\rightarrow \mathbb{R}$. The space of additively separable functions is denoted by
$$\mathcal{D}^\circ_T\triangleq\{H=(H_0,\cdots,H_{T+1})\in \mathcal{D}_T|~H_t~\text{is additively separable for}~ t=0,\cdots,T,~\text{and}~H_{T+1}(\cdot)\equiv0\}.$$
Let $\boldsymbol\mu\triangleq(\boldsymbol\mu_0,\cdots,\boldsymbol\mu_{\tau})$ with $\boldsymbol\mu_t\in \mathbb{R}^{L}_{+}$. We define the operator $\mathcal{L}^{\circ}_T$ on $\mathcal{D}^\circ_T$:
\begin{align}
\mathcal{L}^{\circ}_T H(\textbf{x}_0) \triangleq & \mathbb{E}_0\left[\min_{ \boldsymbol\mu\geq 0}\max_{\textbf{a}\in \mathcal{A}(T)}  \left\{\sum_{t=0}^{T}\big(R_t(\textbf{x}_t,\textbf{a}_t)+ \boldsymbol\mu_t^{\top}(\textbf{b}_t-\textbf{B}_t(\textbf{x}_t,\textbf{a}_t))+ \mathbb{E}[ H_{t+1}(\textbf{x}_{t+1})|\textbf{x}_t,\textbf{a}_t]- H_{t+1}(\textbf{x}_{t+1})\big)\right\}\right]  \label{Finite_Upper_Bound}\\
=&H_0(\textbf{x}_{0})+\mathbb{E}_0\left[\min_{ \boldsymbol\mu\geq 0}\max_{\textbf{a}\in \mathcal{A}(T)} \left\{I_{H}(\textbf{a},\omega,T; \boldsymbol\mu)\right\}\right], \notag
\end{align}
where
$$I_{H}(\textbf{a},\omega,T; \boldsymbol\mu)\triangleq\sum_{t=0}^{T}\big(R_t(\textbf{x}_t,\textbf{a}_t)+  \boldsymbol\mu_t^{\top}(\textbf{b}_t-\textbf{B}_t(\textbf{x}_t,\textbf{a}_t))+ \mathbb{E}[ H_{t+1}(\textbf{x}_{t+1})|\textbf{x}_t,\textbf{a}_t]- H_{t}(\textbf{x}_{t+1})\big).$$

We list the analogous results of Theorem 1, Theorem 2, Theorem 4, and Theorem 5 for finite horizon problem in Theorem 6. Proofs are similar and hence are omitted here.
\begin{thm9}\label{Theorem:FiniteHorzion}
\begin{enumerate}
\item[(a)] (Weak Duality) For any $H\in \mathcal{D}_T$, $ V_0(\textbf{x}_0)\leq \mathcal{L}_T H(\textbf{x}_0)$ for all $\textbf{x}_0\in \mathcal{X}$.

\item[(b)]
(Tighter Bound) For any $H\in\mathcal{D}^\ast_T$, where
 {\small$$\mathcal{D}^\ast_T\triangleq\left\{H\in \mathcal{D}_T: R_t(\textbf{x}_t,\textbf{a}_t)+\beta \mathbb{E}[ H_{t+1}(\textbf{x}_{t+1})|\textbf{x}_t,\textbf{a}_t]\leq H_t(\textbf{x}_{t}) ~~\text{for all}~  \textbf{x}_t\in  \mathcal{X}  ~\text{and}~ \textbf{a}_t\in \bar{\mathcal{A}}(\textbf{x}_t), ~t=0,\cdots,T \right\},$$}then $\max_{\textbf{a}\in \bar{\mathcal{A}}(T)}\{I_H(\textbf{a},\omega,T)\}\leq0$ for every $\omega\in \Omega$; consequently, $V_0(\textbf{x}_0)\leq \mathcal{L}_T H(\textbf{x}_0)\leq H_0(\textbf{x}_0)$ for all $\textbf{x}_0\in \mathcal{X}$.
 \item [(c)] (Strong Duality)  $V_0(\textbf{x}_0)=\mathcal{L}_T V(\textbf{x}_0)$ for all $\textbf{x}_0\in  \mathcal{X}$, where $V=(V_0,\cdots,V_T)$.
%
\item[(d)] (Comparing Lagrangian Bound) For all $\textbf{x}_0\in  \mathcal{X}$, $V_0(\textbf{x}_0)\leq \mathcal{L}_TJ^{\boldsymbol\lambda}((\textbf{x}_0)\leq J^{\boldsymbol\lambda}_0(\textbf{x}_0),$ where $J^{\boldsymbol\lambda}=(J^{\boldsymbol\lambda}_0,\cdots,J^{\boldsymbol\lambda}_T)$.
\item[(e)] (Relaxed Inner Optimization Problem) Suppose that $H\in \mathcal{D}^\circ_T$, i.e., $H_{t}(\textbf{x}_t)=\theta_t+\sum_{n=1}^N H^n_{t}(x_t^n)$,  $\min_{\boldsymbol\mu\geq 0} \max_{\textbf{a}\in \mathcal{A}(T)}I_{J^{\boldsymbol\lambda}}(\textbf{a},\omega,T;\boldsymbol\mu)\leq 0$  for every $\omega\in \Omega$. Consequently, $\mathcal{L}_T J^{\boldsymbol\lambda}(\textbf{x}_0)\leq\mathcal{L}_T^\circ J^{\boldsymbol\lambda}(\textbf{x}_0)\leq J^{\boldsymbol\lambda}_0(\textbf{x}_0)$ for  all $\textbf{x}_0\in  \mathcal{X}$.
\item[(f)] (Duality Gap) Suppose that $H\in \mathcal{D}^\circ_T$, i.e., $H_{t}(\textbf{x}_t)=\theta_t+\sum_{n=1}^N H^n_{t}(x_t^n)$, and Assumptions 1-3 in Section 3.2 hold for every $\omega\in \Omega$. Then for all $\textbf{x}_0\in  \mathcal{X}$,
  \begin{equation}
  \mathcal{L}^\circ_T H(\textbf{x}_0)-\mathcal{L}_TH(\textbf{x}_0)\leq (1+L(T+1)) \max_{n=1,\cdots, N} \sum_{t=0}^T\Gamma^n_t,
  \end{equation}
where
\begin{align*}
\Gamma^n_t=&\sup_{x^n_t\in \mathcal{X}^n, a^n_t \in \mathcal{A}^n(x^n_0)}\left\{ R^n(x^n_t,a^n_t)+\beta \mathbb{E}[ H^n(x^n_{t+1})|x^n_t,a^n_t]- H^n(x^n_{t}) \right\}\\
&- \inf_{x^n_t\in \mathcal{X}^n, a^n_t \in \mathcal{A}^n(x^n_t)}\{ R^n(x^n_t,a^n_t)+\beta \mathbb{E}[ H^n(x^n_{t+1})|x^n_t,a^n_t]- H^n(x^n_{t})\}.
\end{align*}
\end{enumerate}
\end{thm9}

%







\end{document}